\documentclass[11pt]{article}%
\usepackage{amsfonts}
\usepackage{graphicx}
\usepackage{amsmath}
\usepackage{amssymb}
\usepackage{dsfont}
\usepackage{enumitem}
\usepackage[a4paper,bottom=3.5cm,top=2cm,left=2.5cm,right=2.5cm]{geometry}%
\providecommand{\U}[1]{\protect\rule{.1in}{.1in}}
\newcommand\beq{\begin{equation}}
\newcommand\eeq{\end{equation}}
\newcommand{\E}{{\ensuremath{\mathbb E}} }
\newcommand{\bC}{{\ensuremath{\mathbf C}} }
\newcommand{\I}{{\ensuremath{\mathbf 1}} }
\newcommand{\bbX}{{\ensuremath{\mathbb X}} }
\newcommand{\bbZ}{{\ensuremath{\mathbb Z}} }
\newcommand{\bbY}{{\ensuremath{\mathbb Y}} }
\newcommand{\bX}{{\ensuremath{\mathbf X}} }
\newcommand{\bV}{{\ensuremath{\mathbf V}} }
\newcommand{\bU}{{\ensuremath{\mathbf U}} }
\newcommand{\btV}{{\widetilde{\ensuremath{\mathbf V}}} }
\newcommand{\btX}{{\widetilde{\ensuremath{\mathbf X}}} }
\newcommand{\bbtX}{{\widetilde{\ensuremath{\mathbb X}}} }
\newcommand{\btS}{{\widetilde{\ensuremath{\mathbb S}}} }
\newcommand{\bbtS}{{\widetilde{\ensuremath{\mathbb S}}} }
\newcommand{\bS}{{\ensuremath{\mathbb S}} }
\newcommand{\bbS}{{\ensuremath{\mathbb S}} }
\newtheorem{theorem}{Theorem}

\newtheorem{corollary}[theorem]{Corollary}

\newtheorem{lemma}[theorem]{Lemma}

\newtheorem{proposition}[theorem]{Proposition}

\begin{document}

\begin{center}
{\Large \textbf{Bernstein type inequality for a class of dependent random matrices}} \medskip

Marwa Banna$^{a}$, Florence Merlev\`{e}de$^{b}$, Pierre Youssef$^{c}$
\end{center}

$^{a \, b}$ Universit\'{e} Paris Est, LAMA (UMR 8050), UPEMLV, CNRS, UPEC, 5
Boulevard Descartes, 77 454 Marne La Vall\'{e}e, France.

E-mail: marwa.banna@u-pem.fr; florence.merlevede@u-pem.fr

$^{c}$ Department of Mathematical and Statistical sciences. 
University of Alberta, Canada.

Email: pyoussef@ualberta.ca

\bigskip



\begin{center}
\bigskip

\bigskip\textbf{Abstract.} 

\end{center}

 In this paper we obtain a Bernstein type inequality for the sum of  self-adjoint centered  and geometrically absolutely regular random matrices with bounded largest eigenvalue. This inequality can be viewed as an extension to the matrix setting of the Bernstein-type inequality obtained by Merlev\`ede et al. (2009) in the context  of real-valued bounded random variables that are geometrically absolutely regular. The proofs rely on decoupling  the Laplace transform of a sum on a Cantor-like set of random matrices.

\section{Introduction}
\setcounter{equation}{0}

\smallskip
\quad
The analysis of the spectrum of large matrices has known significant development recently due to its important role in several domains. One of the questions is to study the fluctuations of a Hermitian matrix $\bbX$ from its expectation measured by its largest eigenvalue. Matrix concentration inequalities give probabilistic bounds for such fluctuations and provide effective methods for studying several models. The Laplace transform method, which is due to Bernstein in the scalar case, was generalized to the sum of independent Hermitian random matrices by Ahlswede and Winter in \cite{AhWi}. The starting point is that the usual Chernoff bound when we deal with the partial sums associated to real-valued random variables has the following counterpart in the matrix setting: 
\beq \label{LT}
{\mathbb P} \Big (  \lambda_{{\rm max}} \big ( \sum_{i =1 }^n \bbX_i  \big ) \geq x  \Big )  \leq  \inf_{ t >0} \Big \{  {\rm e}^{-t x }  \cdot \E {\rm Tr}\big ({\rm e}^{ t \sum_{i =1}^n \bbX_i  } \big )   \Big \}  
\eeq
(see \cite {AhWi}).  Here and all along the paper, $ (\bbX_i )_{i \geq 1}$ is a family of $d \times d$ self-adjoint random matrices.  
The main problem is then to give a suitable bound for $L_n (t) := \E {\rm Tr}\big ({\rm e}^{ t \sum_{i =1}^n \bbX_i  } \big ) $. In the independent case, starting from  the Golden-Thompson inequality stating that if ${\mathbb A}$ and ${\mathbb B}$ are two self-adjoint matrices,
\[
 {\rm Tr} \big ( {\rm e}^{{\mathbb A}+{\mathbb B}} \big ) \leq   {\rm Tr} \big ( {\rm e}^{{\mathbb A}}  {\rm e}^{{\mathbb B}} \big )  \, ,
\]
Ahlswede  and   Winter have observed that 
\beq \label{iterativeproc}
\E {\rm Tr}\big ({\rm e}^{ t \sum_{i =1}^n \bbX_i  } \big )  \leq  \lambda_{{\rm max}}  ( \E  ({\rm e}^{ t  \bbX_n}  ) ) \cdot \E {\rm Tr}\big ({\rm e}^{ t \sum_{i =1}^{n-1} \bbX_i  } \big )  
\eeq
and gave  a bound for $L_n (t)$ by iterating the procedure above. In \cite{Tr}, Tropp  used Lieb's concavity theorem (see \cite{Li}) 
to improve the bound   on $L_n(t)$ stated in \cite{AhWi}   and obtained Lemma \ref{lmatropp} of Section \ref{sectionpreliminary}. This lemma then allows  to extend to the matrix setting the usual Bernstein inequality for the partial sum associated with independent real-valued random variables.  

Let us mention that an extension of the so-called Hoeffding-Azuma inequality for matrix martingales and of the so-called McDiarmid bounded difference inequality for matrix-valued functions of independent random variables is also stated in \cite{Tr}.

\smallskip

Taking another direction, Mackey et al.  \cite{Ma et al} extended to the matrix setting Chatterjee's technique for developing scalar concentration inequalities  via Stein's method of exchangeable pairs (see \cite{Ch-thesis} and \cite{Ch}), and established Bernstein and Hoeffding inequalities as well as concentration inequalities. Following this approach, Paulin et al. \cite{Pa et al} established a so-called McDiarmid inequality for matrix-valued functions of dependent random variables under conditions on the associated Dobrushin interdependence matrix.

\smallskip

The aim of this paper is to give an extension of the Bernstein deviation inequality when we consider the largest eigenvalue of the partial sums associated with self-adjoint, centered and absolutely regular random matrices with bounded largest eigenvalue. This kind of dependence cannot be compared to the dependence structure imposed in \cite{Ma et al} or in \cite{Pa et al}. 

Note that for dependent random matrices,  the first step given by \eqref{iterativeproc} of the iterative procedure in \cite{AhWi} fails as well as the concave trace function method used in \cite{Tr}. Therefore additional transformations on the Laplace transform have to be made. Even  in the scalar dependent case, obtaining sharp Bernstein-type inequalities is a challenging problem and a dependence structure of the underlying process  has obviously to be precise. For instance, Adamczak \cite{Adam} proved a Bernstein-type inequality for the partial sum associated with  bounded functions of a geometrically ergodic Harris recurrent Markov chain. He showed that even in this context where it is possible to go back to the independent setting by creating random iid cycles, a logarithmic extra factor (compared to the independent case)  cannot be avoided (see Theorem 6 and  Section 3.2 in \cite{Adam}).  

In \cite{MPR1} and \cite{MPR2}, Merlev\`ede et al. considered more general dependence structures than Harris recurrent Markov chains and proved Bernstein-type inequalities for the partial sums associated with bounded real-valued random variables whose  strong mixing coefficients (or $\tau$-dependent coefficients) decrease geometrically or sub-geometrically. Note that in \cite{MPR2}, the case of real-valued random variables that are not necessarily bounded is also treated. The method used in both papers mentioned consists of partitioning the $n$ random variables in blocks indexed by Cantor-type sets
plus a remainder.  The idea is then to control the  log-Laplace transform of each partial sum on the Cantor-type sets. The  log-Laplace transform of the total  partial sum  is then handled with the help of a general result which provides bounds for the
log-Laplace transform of any sum of real-valued random variables (see our Lemma \ref{breta} in the context of random matrices). Obviously, the main step is to obtain a suitable upper bound  
of the  log-Laplace transform of the partial sum on each of  the Cantor-type set. The  dependence structure assumed  in  \cite{MPR1} or  \cite{MPR2} allow the following control: for any index sets $Q$ and $Q'$ of natural numbers such that  
$Q \subset [1, p]$ and $Q' \subset [n+p, \infty)$ and any $t >0$,
\beq \label{deco}
\E \big ( {\rm e}^{t \sum_{i \in Q} X_i}  {\rm e}^{t \sum_{i \in Q'} X_i}  \big )  \leq  \E \big ( {\rm e}^{t \sum_{i \in Q} X_i}  \big ) \E \big (  {\rm e}^{t \sum_{i \in Q'} X_i}  \big )  + \varepsilon (n)  \Vert   {\rm e}^{t \sum_{i \in Q} X_i} \Vert_{\infty}  \Vert   {\rm e}^{t \sum_{i \in Q'} X_i} \Vert_{\infty} \, ,
\eeq
where $\varepsilon(n)$ is a sequence of positive real numbers depending on the dependent coefficients.  The binary tree structure of the Cantor-type  sets allows to iterate the decorrelation procedure above-mentioned allowing then to suitably handle  the log-Laplace transform of the partial sum on each of  the Cantor-type set.

In the random matrix setting, iterating a procedure as \eqref{deco} cannot lead to suitable exponential inequalities essentially due to the fact that the extension of the Golden-Thompson inequality to three or more Hermitian matrices fails, and then the extension of the exponential inequalities stated in  \cite{MPR1} and \cite{MPR2} to the matrix setting is not straight forward. To benefit of the ideas developed in \cite{AhWi} or in \cite{Tr}, we shall rather bound the  log-Laplace transform of the partial sum indexed by  a Cantor-type set, say $K$, by the  log-Laplace transform of the sum of $2^\ell $  independent and self-adjoint random matrices   plus a small error term (here $\ell$ depends on the cardinality of $K$).  Lemma \ref{lemme1} is in this direction and  can be viewed as a decoupling lemma for the Laplace transform in the matrix setting.  As we shall see, a well-adapted dependence structure allowing  such a procedure is the absolute regularity structure. Indeed, the Berbee's coupling lemma (see Lemma \ref{berbeelma} below)  allows a "good coupling" in terms of absolute regular coefficients (see the definition \eqref{defbeta}) even when the underlying random variables take values in a high dimensional space (working with $d \times d $  random matrices can be viewed as working with random vectors of dimension $d^2$).  The decoupling lemma \ref{lemme1} associated with additional coupling arguments   will then allow us to prove our key Proposition \ref{propinter1} giving a bound for the Laplace transform of the partial sum indexed by Cantor-type set of self-adjoint random matrices. As we shall see, our method allows to extend the scalar Bernstein type inequality given in \cite{MPR1} to the matrix setting. 

Our paper is organized as follows. In Section \ref{sectionmain}, we introduce some notations and definitions and state our Bernstein-type inequality  for the class of random matrices we consider (see  Theorem \ref{mainthm}). Section \ref{sectionappli} is devoted to some examples of matrix models where this  Bernstein-type inequality applies. The proof of the main result is given in Section \ref{section proof}.

\section{Main Result}    \label{sectionmain}
For any   $d \times d$ matrix $\bbX = [ (\bbX )_{i,j}]_{i,j=1}^d$ whose entries belong to  ${\mathbb K} = {\mathbb R}$ or ${\mathbb C}$, we associate its corresponding vector  $\bX$ in ${\mathbb K}^{d^2}$ whose coordinates are the entries of $\bbX$ i.e.
\[
\bX = \big ( (\bbX )_{i,j} \, ,  1 \leq  i \leq  d \big )_{1 \leq  j \leq  d} \, .
\]
Therefore $\bX = \big ( X_i \, , 1 \leq  i \leq  d^2 \big )$ where
\[
X_i =  (\bbX )_{i-(j-1)d, j } \ \text{ for $(j-1) d+1 \leq  i \leq  jd$}\, ,
\]
and $\bX$ will be called the vector associated with $\bbX$. Reciprocally, given   $\bX= \big ( X_\ell \, , 1 \leq  \ell \leq  d^2 \big )$ in ${\mathbb K}^{d^2}$ we shall associate a $d \times d$ matrix  $\bbX$ by setting
\[
\bbX = \big [ (\bbX)_{i,j} \big ]_{i,j =1}^n \text{ where } (\bbX)_{i,j} = X_{i+(j-1)d} \, .
\]
The matrix $\bbX$ will be referred to as the matrix associated with $\bX$. 
\medskip

In all the paper we consider a family $(\bbX_i)_{i \geqslant 1}$ of $d \times d$ self-adjoint random matrices  whose entries are defined on a probability space $(\Omega, {\mathcal A}, {\mathbb P})$ and with values in ${\mathbb K}$,  and that are geometrically absolutely regular in the following sense. Let 
\beq \label{defbeta}
\beta_0 = 1 \ \text{ and } \ \beta_k = \sup_{j \geqslant 1} \beta ( \sigma ( \bX_i \, , \, i \leq  j ),  \sigma ( \bX_i \, , \, i \geqslant j +k) )  \, , \text{ for any $k \geq 1$,}
\eeq
where \[
\beta ({\mathcal{A}},{\mathcal{B}})=\frac{1}{2}\sup \big \{\sum_{i\in
I}\sum_{j\in J}|\mathbb{P}(A_{i}\cap B_{j})-\mathbb{P}(A_{i})\mathbb{P} (B_{j})|\big \} \, , 
\]
 the maximum being taken over all finite partitions $%
(A_{i})_{i\in I}$ and $(B_{i})_{i\in J}$ of $\Omega $ respectively with
elements in ${\mathcal{A}}$ and ${\mathcal{B}}$.  

The $(\beta_k)_{k \geqslant 0}$ are usually called the coefficients of absolute regularity of the sequence of vectors  $(\bX_i)_{i \geqslant 1}$ and we shall assume in this paper that they decrease geometrically in the sense that  there exists $c >0$ such that for any integer $k \geqslant 1$, 
\beq \label{condbeta}
\beta_k = \sup_{j \geqslant 1} \beta ( \sigma ( \bX_i \, , \, i \leq  j ),  \sigma ( \bX_i \, , \, i \geqslant j +k) )  \leq  {\rm e}^{-c (k-1) } \, ,
\eeq
Note that  the $\beta_k$ coefficients have been introduced by Kolmogorov and Rozanov \cite{KR}  and  even if they are more restrictive than the so-called Rosenblatt strong mixing coefficients $\alpha_k$ they can be computed in many situations. For instance, we refer  to the work by Doob \cite{Do} for sufficient conditions on Markov chains to be geometrically absolutely regular or by Mokkadem \cite{Mo} for mild conditions ensuring vector ARMA processes to be also geometrically $\beta$-mixing.
\medskip

 In all the paper, we will assume that the underlying probability space $(\Omega,{\cal A}, {\mathbb P} )$  is rich enough  to contain a sequence   
$(\epsilon_i)_{i \in {\mathbb Z}} = (\delta_i, \eta_i)_{i \in {\mathbb Z}}$ of iid random variables with uniform distribution over $[0,1]^2$, independent of $(\bX_i)_{i \geqslant 0}$. In addition, the following notations will be used : $\log x := \ln x$,  $\log_2 x = \frac{\log x}{ \log 2}$, we write ${\mathbf  0}$ for the zero matrix and ${\mathbb  I}_d$ for the $d \times d$ identity matrix, we use the curly inequalities  to denote the semidefinite ordering i.e.  ${\bf 0  \preceq \bbX }$ means that $\bbX$ is positive semidefinite.

\smallskip

\begin{theorem} \label{mainthm} Let $(\bbX_{i})_{i \geqslant 1}$ be a family of self-adjoint random matrices of size $d$. Assume  that  \eqref{condbeta} holds and that there exists a positive constant $M$ such that for any $i \geq 1$, 
\beq \label{condmajlamdamax}
\E(\bbX_i) = {\bf 0} \quad \text{and} \quad \lambda_{{\rm max}} ( \bbX_i ) \leq  M \quad \text{almost surely.}
\eeq
Then there exists a universal positive constant $C$ such that for any $x>0$ and any integer $n \geqslant 2$,
\[
{\mathbb P} \Big (  \lambda_{{\rm max}} \big ( \sum_{i =1 }^n \bbX_i  \big ) \geq x  \Big )  \leq  d \exp \Big (-\frac{C x^{2}}{   v^{2}n+  c^{-1}M^{2} +  xM \gamma (c, n)} \Big )\,.  
\]
where
\begin{equation}
v^{2}= \sup_{K\subseteq \{1, \dots, n\} }\frac{1}{\mathrm{Card}K}%
 \lambda_{{\rm max}}  \Big ( \E \big ( \sum_{i \in K} \bbX_i  \big )^2 \Big )  \label{hypovarcont}
\end{equation}%
and 
\begin{equation}
 \gamma (c, n) =  \frac{\log n}{\log 2} \max \Big (   2 , \frac{32 \log n}{c  \log 2} \Big ) \, .
 \end{equation}\label{gamma}
\end{theorem}
In the definition of $v^2$ above, the maximum is taken over all nonempty  subsets $K\subseteq \{1, \dots, n\}$.

\medskip
To prove the deviation inequality stated in Theorem \ref{mainthm}, we shall use the matrix Chernoff  bound \eqref{LT}. The theorem will then follow from the following control of the matrix log-Laplace transform that is proved in Section \ref{sectionend}:  Under the conditions of Theorem \ref{mainthm},  for any positive $t$ such that $\displaystyle tM < 1/\gamma (c, n) $, we have
\[
\log \E   {\rm Tr} \Big (\exp \big (  t \sum_{i =1}^n\bbX_i    \big ) \Big ) \leq  \log d +  \frac{     t^2 n  \big ( 15 v + 2  M/ (cn)^{1/2} \big )^2 }{1 - t M \gamma (c, n)  } \, .
\]
As proved in Section 4.2.4 of \cite{Ma et al}, this inequality together with Jensen's inequality leads to the following upper bound for the expectation of the largest eigenvalue of $\sum_{i =1}^n\bbX_i $: Under the conditions of Theorem \ref{mainthm},
\[
\E  \lambda_{{\rm max}} \Big ( \sum_{i =1 }^n \bbX_i  \Big ) \leq 30v \sqrt{n \log d}  + 4 M c^{-1/2}  \sqrt{ \log d}  + M \gamma(c, n) \log d  \, .
\]

\section{Applications} \label{sectionappli}

Let $(\tau_k)_k $ be a stationary sequence of real-valued random variables such that  $\| \tau_1 \|_{\infty} \leq 1$ a.s.  Consider a family $(\bbY_k)_k$ of independent real and symmetric $d \times d$ random matrices which is independent of $(\tau_k)_k$. 
For any $i=1, \ldots,n$, let $\bbX_i= \tau_i \bbY_i$ and note that in this case
$$
\beta_k= \beta ( \sigma ( \tau_i \, , \, i \leq 0 ),  \sigma ( \tau_i \, , \, i \geq k) )  \, .
$$
\begin{corollary}
Assume that there exists a positive constant $c$ such that $\beta_k\leq e^{-c(k-1)}$  for any $k\geq 1$
and suppose that each random matrix $\bbY_k$ satisfies
$$
\E\bbY_k={\bf 0} \quad \text{,} \quad  \lambda_{{\rm max}} (\bbY_k) \leq M \quad \text{and} \quad  \lambda_{{\rm min}} (\bbY_k) \geq - M \quad \text{almost surely.}
$$
Then for any $t>0$ and any integer $n \geq 2$, 
\[
{\mathbb P} \Big (  \lambda_{{\rm max}} \big ( \sum_{k =1 }^n \tau_k \bbY_k  \big ) \geq t  \Big )  
\leq d \exp \Big (-\frac{C t^{2}}{ nM^2 \E(\tau_0^2) + M^{2} +t M (\log n)^2} \Big )\,  ,
\]
where $C$ is a positive constant depending only on $c$. 
\end{corollary}

\noindent
\textit{Proof.} The above corollary follows by noting that for any $K\subseteq \{1, \dots, n\}$
$$
\Sigma_K:=\E \big ( \sum_{k \in K} \tau_k \bbY_k  \big )^2 = \sum_{k \in K} \E(\tau_k^2) \E(\bbY_k^2)=  \E(\tau_0^2) \sum_{k \in K}  \E(\bbY_k^2),
$$
which, by Weyl's inequality, implies that $\lambda_{{\rm max}}(\Sigma_K) \leq M^2 \text{Card}(K) \E(\tau_0^2)$. Therefore, we infer that $v^2 \leq M^2 \E(\tau_0^2)$. 
$\square $

\bigskip

We consider now another model for which Theorem \ref{mainthm} can be applied. Let $(X_k)_{k \in {\mathbb Z} }$ be a  geometrically absolutely regular sequence of real-valued centered random variables. That is, there exists a positive constant $c_{0}$ such that for any $k\geq 1$,
\begin{equation}\label{mixing-rate-example}
 \sup_{\ell \in {\mathbb Z}} \beta\big(\sigma ( X_i \, , \, i \leq \ell ),  \sigma ( X_i \, , \, i \geq k + \ell) \big)  \leq e^{-c_0 (k-1)} \, .
\end{equation}
For any $i=1, \ldots ,n$,  let $\bbX_i$ be the $d \times d$ random matrix defined by $\bbX_i =\bC_i \bC_i^T - \E (\bC_i \bC_i^T)$ where $ \bC_i=( X_{(i-1)d+1}, \ldots , X_{id})^T $. 
Note that 
in this case,  
$$
\beta_k =  \sup_{\ell \in {\mathbb Z}}  \beta \big( \sigma ( \bC_i , i \leq \ell),  \sigma ( \bC_i \, , \, i \geq \ell +k) \big) 
\leq e^{-c_0 d(k-1)}  \, .
$$
for any $k\geq 1$.
\begin{corollary}
Assume that $(X_k)_k$ satisfies \eqref{mixing-rate-example}. Suppose in addition that there exists a positive constant $M$  satisfying $ \sup_k \Vert X_k \Vert_{\infty} \leq M$ a.s. Then,  for any  $x >0$ and any integer $n \geq 2$ 
\[
{\mathbb P} \Big (  \lambda_{{\rm max}} \big ( \sum_{i =1 }^n \bbX_i  \big ) \geq x  \Big )  \leq d \exp \Big (-\frac{C x^{2}}{  n d M^4 + dM^4 +x M^2  ( d  \log n + \log^2n ) } \Big )\,,  
\]
where $C$ is a positive constant depending only on $c_0$. 
\end{corollary}
\textit{Proof.} For any $i\in \{1, \ldots, n\}$,  note that $\lambda_{{\rm max}}(\bbX_i) \leq \lambda_{{\rm max}} (\bC_i \bC_i^T)$ implying that  $ \lambda_{{\rm max}}(\bbX_i) \leq  d M^2 $ a.s. 
To get the desired result, it remains to control $v^2$. We have for any $K \subseteq \{ 1, \ldots , N\}$,
\[
\Sigma_K:=\E \big( \sum_{i\in K} \bbX_i\big)^2 
= 
\sum_{i,j \in K} \text{Cov} (\bC_i\bC_i^T, \bC_j\bC_j^T)
\]
and we  note that the $(k,\ell)^{th}$ component of $\Sigma_K$ is 
\[
(\Sigma_K)_{k,\ell}=\Big[ \E \big( \sum_{i\in K} \bbX_i\big)^2 \Big]_{k,\ell} 
= \sum_{i,j \in K} \sum_{s=1}^d \text{Cov} \big( X_{(i-1)d+k}\, X_{(i-1)d+s}\, , X_{(j-1)d+s} \,X_{(j-1)d+\ell}\big) \, .
\]
Therefore we infer by Gerschgorin's theorem  that 
\begin{multline*}
\big| \lambda_{{\rm max}} \big (\Sigma_K \big) \big| \\
\leqslant \sup_k \sum_{\ell=1}^d |(\Sigma_K)_{k,\ell}|
\leqslant \sup_k \sum_{i,j \in K} \sum_{\ell=1}^d \sum_{s=1}^d \big| \text{Cov} \big( X_{(i-1)d+k}\, X_{(i-1)d+s}\, , X_{(j-1)d+s} \,X_{(j-1)d+\ell}\big)\big| \, .
\end{multline*}
After tedious computations involving Ibragimov's covariance inequality (see \cite{Ibra}), we infer that $v^2 \leq c_1d M^4$ where $c_1$ is a  positive constant depending only on $c_0$. Applying Theorem \ref{mainthm} with these upper bounds ends the proof.
$\square $

\section{Proof of  Theorem \ref{mainthm}} \label{section proof}

The proof of  Theorem \ref{mainthm} being very technical, it is divided into several steps. In Section \ref{sectionpreliminary}, we first collect some technical preliminary lemmas that will be necessary all along the proof. In  Section \ref{section key}, we give the main ingredient to prove our Bernstein-type inequality, namely: a bound for the Laplace transform of the partial sum, indexed by a suitable Cantor-type set, of the self-adjoint random matrices under consideration (see  Proposition \ref{propinter1} and  Section \ref{sectionCantor} for the construction of this suitable Cantor-set). As quoted in the introduction, this key result is based on a decoupling lemma which is stated in Section \ref{sectionfundamentallemma}. The proof of Theorem  \ref{mainthm} is completed in Section \ref{sectionend}.

\subsection{Preliminary materials} \label{sectionpreliminary}

The following lemma is due to Tropp \cite{Tr}.  Under the form stated below, it is a combination of his Lemmas 3.4 and 6.7 together with the proof of his Corollary 3.7.

\begin{lemma} \label{lmatropp} 
Let $K$ be a finite subset of positive integers. Consider a family $({\mathbb U}_k)_{k \in K}$ of $d \times d$ self-adjoint random matrices  that are mutually independent. Assume that for any $k \in K$,
\[
\E ( {\mathbb U}_k) = {\mathbf 0} 
\quad  \text{and} \quad
\lambda_{{\rm max}} ( {\mathbb U}_k)  \leq  B  \  \text{a.s.}
\]
where $B$ is a positive constant. Then for any $t>0$,
\begin{equation}  \label{psi2}
\E {\rm Tr}\big ({\rm e}^{t  \sum_{k \in K}{\mathbb U}_k  } \big ) 
 \leq  d  \exp  \Big (    t^2 g(t B) \lambda_{{\rm max}}  \Big ( \sum_{k \in K} \mathbb{E} ({\mathbb U}_k^2) \Big )  \Big ) \, ,
\end{equation}
where $g(x) =x^{-2}( {\rm e}^{x} -x-1)$. 
\end{lemma}
The next lemma is an adaptation of Lemma 3 in \cite{MPR2} to the case of the log-Laplace transform  of any sum of $d \times d$ self-adjoint random matrices.  
\begin{lemma}
\label{breta}
Let ${\mathbb U}_0, {\mathbb U}_1, \ldots$ be a sequence of $d \times d$ self-adjoint random matrices. Assume that there exists positive constants $\sigma _{0},\sigma
_{1},\ldots $ and $\kappa_{0},\kappa_{1},\ldots $ such that, for any $i\geq 0$ and
any $t$ in $[0,1/\kappa_{i}[$,
\begin{equation*}
\log \E {\rm Tr}\big ({\rm e}^{t \, {\mathbb U}_i  } \big ) \leq  C_d +  (\sigma _{i}t)^{2}/(1-\kappa_{i}t)\,,
\end{equation*}
where $C_d$ is a positive constant depending only on $d$.  Then, for any positive $n$ and any $t$ in $[0, 1/(\kappa_0 + \kappa_1 + \cdots + \kappa_n)[$%
, 
\begin{equation*}
\log\E {\rm Tr}\big ({\rm e}^{ t \sum_{k =0}^n{\mathbb U}_k  } \big )  \leq   C_d +  (\sigma t)^{2}/(1-\kappa t),
\end{equation*}%
where $\sigma =\sigma _{0}+\sigma _{1}+\cdots +\sigma _{n}$ and $%
\kappa = \kappa_{0}+ \kappa_{1}+\cdots + \kappa_{n}$.
\end{lemma}
\noindent{\bf Proof.} Lemma \ref{breta} follows
from the case $n=1$ by induction on $n$. For any $t \geq 0$, let
\[
L(t) = \log \E {\rm Tr}\big ({\rm e}^{t  ( {\mathbb U}_0 + {\mathbb U}_1 )  } \big ) 
\]
and notice that by the Golden-Thompson inequality,
\beq \label{GTL(t)}
L(t) \leq  \log \E {\rm Tr}\big ({\rm e}^{t  \,  {\mathbb U}_0 } {\rm e}^{t  \,   {\mathbb U}_1  } \big )  \, .
\eeq
Define  the functions $\gamma _{i}$ by 
\begin{equation*}
\gamma _{i}(t)=(\sigma _{i}t)^{2}/(1-\kappa_{i}t)\ \hbox{ for }t\in \lbrack
0,1/\kappa_{i}[\ \hbox{ and }\gamma _{i}(t)=+\infty \ \hbox{ for }t\geq 1/\gamma_{i} \, ,
\end{equation*}
and recall the non-commutative H\"older inequality (see for instance exercise 1.3.9 in \cite{Tao12}): if $A$ and $B$ are $d \times d$ self-adjoint random matrices then, for any $1 \leq p,q \leq \infty$ with $p^{-1}+q^{-1} = 1$,
\beq \label{ineSchatten}
\vert {\rm Tr} ( A B) \vert  \leq \Vert A \Vert_{{\mathcal S}^p} \Vert B \Vert_{{\mathcal S}^q}  \, ,
\eeq
where $\Vert A \Vert_{{\mathcal S}^p}  = \Vert ( \lambda_i(A) )_{i=1}^d \Vert_{\ell_n^p}=\Big (  \sum_{i=1}^d \vert \lambda_i(A)\vert^p \Big )^{1/p} $ (resp. $\Vert B \Vert_{{\mathcal S}^q}$) is the $p$-Schatten norm of $A$ (resp the $q$-Schatten norm of $B$). 

Starting from  \eqref{GTL(t)} and applying \eqref{ineSchatten} with  $A= {\rm e}^{t  \,  {\mathbb U}_0 }$ and $B= {\rm e}^{t  \,  {\mathbb U}_0 }$, we derive that for any $t >0$ and any $p \in ]1, \infty[$
\[
L(t) \leq  \log  \E  \Big ( \Vert {\rm e}^{t  \,  {\mathbb U}_0 } \Vert_{{\mathcal S}^{p}}  \Vert {\rm e}^{t  \,  {\mathbb U}_1 } \Vert_{{\mathcal S}^{q}}  \Big ) \, ,
\]
which gives by applying H\"older's inequality
\[
L(t) \leq p^{-1} \log  \E \Vert {\rm e}^{t  \,  {\mathbb U}_0 } \Vert^{p}_{{\mathcal S}^{p}} + q^{-1} \log \E  \Vert {\rm e}^{t  \,  {\mathbb U}_1} \Vert^{q}_{{\mathcal S}^{q}} \, .
\]
Observe now that since ${\mathbb U}_0$ is self-adjoint
\[
\Vert {\rm e}^{t  \,  {\mathbb U}_0 } \Vert^{p}_{{\mathcal S}^{p}} = \sum_{i=1}^d \vert \lambda_i( {\rm e}^{t  \,  {\mathbb U}_0 } )\vert^p  
=  \sum_{i=1}^d \lambda_i( {\rm e}^{t  p \,  {\mathbb U}_0 } ) =  {\rm Tr}\big ({\rm e}^{t  p  \, {\mathbb U}_0  } \big ) \, , \]
and similarly $\Vert {\rm e}^{t  \,  {\mathbb U}_1 } \Vert^{q}_{{\mathcal S}^{q}}  =  {\rm Tr}\big ({\rm e}^{t  q  \,  {\mathbb U}_1  } \big )$.  So, overall, 
\beq \label{ine2L(t)}
L(t) \leq p^{-1} \log  \E  {\rm Tr}\big ({\rm e}^{t  p  \, {\mathbb U}_0  } \big )+ q^{-1} \log \E {\rm Tr}\big ({\rm e}^{t  q  \,  {\mathbb U}_1  } \big ) \, .
\eeq
For any  $t$ in $[0,1/\kappa[$, take $u_t=(\sigma _{0}/\sigma )(1-\kappa t)+\kappa_{0}t$
(here $\kappa=\kappa_{0}+\kappa_{1}$ and $\sigma =\sigma _{0}+\sigma _{1}$). With this
choice $1-u_t=(\sigma _{1}/\sigma )(1-\kappa t)+\kappa_{1}t$, so that $u_t$ belongs to $%
]0,1[$. Applying Inequality \eqref{ine2L(t)} with $p = 1/u_t$,  we get that for any  $t$ in $[0,1/\kappa[$,
\[
L(t) \leq u_t \gamma_0 (t /u_t) + (1-u_t) \gamma_1 (t /(1-u_t))=(\sigma t)^{2}/(1-\kappa t) \, , 
\]
which completes
the proof of Lemma \ref{breta}. $\square $

\medskip

Next lemma allows coupling and is due to Berbee \cite{Be}.
\begin{lemma} \label{berbeelma}
Let $X$ and $Y$ be two random variables defined on a probability space $( \Omega, {\mathcal A}, {\mathbb P})$ and taking their values in Borel spaces $B_1$ and $B_2$ respectively. Assume that 
$( \Omega, {\mathcal A}, {\mathbb P})$ is rich enough to contain a random variable $\delta$ with uniform distribution over $[0,1]$ independent of $(X,Y)$. Then there exists a random variable $Y^*=f(X,Y,\delta)$ where $f$ is a measurable function from $B_1 \times B_2 \times [0,1]$ into $B_2$ such that $Y^*$ is independent of $X$, has the same distribution as $Y$ and
\[
{\mathbb P} ( Y \neq Y^*) = \beta (\sigma(X), \sigma(Y)) \, .
\]
\end{lemma}
Let us note that the $\beta$-mixing coefficient $\beta (\sigma(X), \sigma(Y))$ has the following equivalent definition:
\beq \label{equivalentdefbeta}
\beta (\sigma(X), \sigma(Y)) = \frac{1}{2} \Vert P_{X,Y} - P_X \otimes P_Y  \Vert \, ,
\eeq
where $P_{X,Y} $ is the joint distribution of $(X,Y)$ and  $P_X $ and $P_Y$ are respectively the distributions of $X$ and $Y$ and,  for two positive measures $\mu$ and $\nu$, the notation  $\Vert \mu - \nu \Vert$ denotes the total variation of 
$\mu - \nu$.

\subsection{A key result} \label{section key}

The next proposition is the main ingredient to prove  Theorem \ref{mainthm}. It is based on a suitable construction of a subset $K_{A}$ of $\{1,\dots ,A\}$ for which it is possible to give a good upper bound for the Laplace transform of $\sum_{i \in K_{A}}\bbX_i  $. Its proof is based on the decoupling Lemma \ref{lemme1} below that allows  to compare $ \E {\rm Tr}\big ({\rm e}^{t \sum_{i \in K_{A}}\bbX_i    } \big )$ with the same quantity but replacing $\sum_{i \in K_{A}}\bbX_i  $ with a sum of independent blocks.

\begin{proposition} \label{propinter1} Let $(\bbX_{i})_{i \geqslant 1}$ be as in Theorem \ref{mainthm}. Let $A $ be a positive integer larger than $2$. Then there exists a subset $K_{A}$ of $\{1,\dots ,A\}$ with $\mathrm{Card}(K_{A})\geqslant A/2$, such
that for any positive $t$ such that $t M \leq   \min  \big ( \frac{1}{2} , \frac{c \log 2}{32 \log A} \big ) $, 
\beq \label{resultpropinter1}
\log \E {\rm Tr}\Big ({\rm e}^{t  \sum_{i \in K_{A}}\bbX_i   } \Big )   \leq  \log d +  4 \times 3.1 t^2 A v^2  
 +   \frac{9 (tM)^2}{ c }   {\rm e}^{ -3c/(32 tM) }     \,,
\eeq
where $v^2$ is defined in \eqref{hypovarcont}. 
\end{proposition}
The proof of this proposition is  divided into several steps.

\subsubsection{Construction of a Cantor-like subset $K_A$} \label{sectionCantor}
 As in \cite{MPR1} and \cite{MPR2}, the set $K_{A}$ will be a finite
union of $2^{\ell }$ disjoint sets of consecutive integers with same
cardinality spaced according to a recursive `Cantor"-like construction.  Let 
\[
\delta = \frac{\log 2}{2 \log A} \ \text{ and } \ \ell:=\ell_A = \sup \{  k \in {\mathbb N}^* \, : \, \frac{A \delta (1 - \delta)^{k-1}}{2^{k}} \geqslant 2\} \, .
\]
Note that $\ell \leq  \log A/ \log 2$ and $\delta \leq  1/2$ (since $A \geqslant 2$).  Let $n_0 = A$ and for any  $j \in \{ 1, \dots,  \ell \}$, 
\beq \label{defdjnj}
n_j = \big \lceil  \frac{A (1 - \delta)^{j}}{2^{j}}  \big \rceil  \ \text{ and } \ d_{j-1} = n_{j-1} - 2 n_j \, . 
\eeq
For any nonnegative $x$, the notation $\lceil x \rceil $ means the smallest integer which is larger or equal to $x$.  Note that for any $j \in \{ 0, \dots,  \ell -1 \}$,
\beq \label{restdj}
d_{j} \geqslant \frac{A \delta (1 - \delta)^{j}}{2^{j}}  -2 \geqslant \frac{A \delta (1 - \delta)^{j}}{2^{j+1}}  \, ,
\eeq
where the last inequality comes from the definition of $\ell$. Moreover,
\beq \label{restnj}
n_\ell \leq  \frac{A (1 - \delta)^{\ell}}{2^{\ell}}  +1 \leq  \frac{A  (1 - \delta)^{\ell}}{2^{\ell-1}}  \, ,
\eeq
where the last inequality comes from the fact that $\displaystyle  \frac{A  \delta(1 - \delta)^{\ell-1}}{2^{\ell}} \times \frac{1-\delta}{\delta} \geqslant 2$ by the definition of $\ell$ and the fact that $\delta \leq  1/2$

To construct $K_A$ we proceed as follows. At the first step, we divide the set $\{1, \dots ,A\}$ into 
three disjoint subsets of consecutive integers: $I_{1,1}$, $I_{0,1}^*$ and $I_{1,2}$. These subsets are such that ${\rm Card}(I_{1,1})={\rm Card}(I_{1,2})=n_1$ and ${\rm Card}(I_{0,1}^*)=d_0$.  At the second step, each of the sets of integers 
$I_{1,i}$, $ i=1,2$, is divided into three disjoint subsets of consecutive integers as follows: for any $i=1,2$,  $I_{1,i} = I_{2, 2i-1} \cup I_{1,i}^* \cup I_{2, 2i}$ where ${\rm Card}(I_{2,2i-1})={\rm Card}(I_{2,2i})=n_2$ and ${\rm Card} ( I_{1,i}^* )=d_1$. 
Iterating this procedure we have constructed after $j$ steps ($1 \leq  j\leq
\ell_{A}$),  $2^{j}$ sets of consecutive integers, $I_{j, i}$, $i=1, \dots, 2^j$, each of cardinality $n_j$ such that $a_{j,2k} - b_{j,2k-1} - 1=d_{j-1}$ for any $k = 1, \dots, 2^{j-1}$, where $a_{j,i} = \min \{ k \in I_{j,i}\}$ and $b_{j,i} = \max \{ k \in I_{j,i}\}$. 
Moreover if, for any $i =1, \dots, 2^{j-1}$, we set $I_{j-1,i}^* = \{ b_{j,2i-1}  +1, \dots, a_{j,2i} -1 \}$, then $I_{j-1, i}=I_{j, 2i-1}\cup I_{j-1,i}^* \cup I_{j,2 i} $. After $\ell$ steps we then have constructed 
$2^{\ell}$ sets of consecutive integers, $I_{\ell, i}$, $i=1, \dots, 2^{\ell}$, each of cardinality $n_{\ell}$ such that $I_{\ell, 2i-1}$ and $I_{\ell, 2i}$ are spaced by $d_{\ell -1}$ integers. 
The set of consecutive integers $K_A$ is then defined by
\begin{equation*}
K_{A}=\bigcup_{k=1}^{2^{\ell }}I_{\ell ,k}\, . 
\end{equation*}
Note that 
\[
\{1, \dots , A \} = K_{A} \cup ( \cup_{j=0}^{\ell -1} \cup_{i=1}^{2^j} { I}^*_{j,i}) 
\]
Therefore
\[
{\rm Card} ( \{1, \dots , A \}  \setminus  K_{A} ) =  \sum_{j=0}^{\ell -1} \sum_{i=1}^{2^j} {\rm Card} ( { I}^*_{j,i} ) = \sum_{j=0}^{\ell -1} 2^j d_j  = A - 2^{\ell} n_{\ell}\, .
\]
But
\beq \label{evident}
A - 2^{\ell} n_{\ell} \leq  A  \big ( 1 - (1- \delta)^{\ell}\big ) = A \delta \sum_{j=0}^{\ell -1} (1 - \delta)^j \leq  A \delta  \ell \leq  \frac{A}{2} \, .
\eeq
Therefore $A \geqslant {\rm Card} (  K_{A} ) \geqslant A/2$. 

In the rest of the proof, the following notation will be also useful: for any $k \in \{0, 1 , \dots, {\ell}\}$ and  any $j \in \{1, \dots, 2^{k}\}$, let
\beq \label{not1}
K_{k,j} := K_{A,k,j}= \bigcup_{ i=(j-1)2^{\ell - k} +1 }^{j 2^{\ell -k}} I_{\ell , i}
\eeq
Therefore $K_{0,1}= K_{A}$ and, for any  $j \in \{1, \dots, 2^{\ell}\}$, $K_{\ell,j} =I_{\ell , j}$. Moreover,  for any  $k \in \{ 1 , \dots, {\ell}\}$ and  any $j \in \{1, \dots, 2^{k-1}\}$, there are exactly $d_{k-1}$ integers between $K_{k,2j-1} $ and $K_{k,2j} $.

\subsubsection{A fundamental decoupling lemma} \label{sectionfundamentallemma}

We start by introducing some notations,  then we state the decoupling  Lemma \ref{lemme1} below that is fundamental to prove Proposition \ref{propinter1}. Let $K_A$ be defined as in Step 1. In the rest of the proof, we will adopt the following notation. For any  integer $m \in \{0, \dots, \ell\}$,  $(\bV_j^{(m)})_{1 \leq  j \leq  2^m}$ will denote  a family of  
$2^m$ mutually independent random vectors defined on $({\Omega, {\mathcal A}, {\mathbb P}})$, each of dimension $s_{d, \ell, m}:=d^2\mathrm{Card}( K_{m,j})= d^2 2^{\ell -m } n_\ell$ and such that 
\beq \label{defofVjm}
\bV_j^{(m)} =^{\mathcal D} (\bX_i \, , i \in K_{m,j} ) \, .
\eeq
The existence of such a family is ensured by the Skorohod lemma (see \cite{Sk}). Indeed since $( \Omega, {\mathcal A}, {\mathbb P})$  is assumed to be large enough to contain a sequence $(\delta_i)_{i \in {\mathbb Z}}$ of iid random variables uniformly distributed on $[0,1]$ and independent of the sequence $(\bX_i)_{ i \geqslant 0}$, there exist measurable functions 
$f_j $ such that the vectors $\bV_j^{(m)}  = f_j \big (  (\bX_i \, , i \in K_{m,k})_{k=1, \ldots, j}  , \delta_j  \big )$, $j=1, \ldots, 2^m$,   are independent and satisfy \eqref{defofVjm}.

Let  $\pi_i^{(m)}$ be the $i$-th canonical projection from ${\mathbb K}^{s_{d, \ell, m}}$ onto ${\mathbb K}^{d^2}$, namely: for any vector ${\mathbf x} = ({\mathbf x}_i \, , i \in K_{m,j} ) $ of ${\mathbb K}^{s_{d, \ell, m}}$,  $\pi_i^{(m)} ( {\mathbf x} ) = {\mathbf x}_i$. 
For any 
$i \in K_{m,j} $, let 
\beq \label{1rstdef}
\bX_{j}^{(m)} (i) = \pi_i^{(m)} ( \bV_j^{(m)} ) \ \text{ and } \   {\mathbb S}_{j}^{(m)} = \sum_{i \in K_{m,j}} \bbX_{j}^{(m)} (i) \, ,
\eeq
where $\bbX_{j}^{(m)} (i)$ is the $d \times d$ random  matrix  associated with $\bX_{j}^{(m)} (i)$ (recall that this means that the $(k, \ell)$-th entry of $\bbX_{j}^{(m)} (i)$ is the 
$( (\ell-1)d +k)$-th coordinate  of the vector $\bX_{j}^{(m)} (i)$).

With the above notations, we have 
\beq \label{trivial} 
\E {\rm Tr}\big ({\rm e}^{t \sum_{i \in K_{A}}\bbX_i    } \big ) =  \E {\rm Tr}\big ({\rm e}^{t   \, \bbS_{1}^{(0)}  } \big )  \, .
\eeq 

\medskip

We are now in position to state the following lemma which will be a key step in the proof of Proposition \ref{propinter1} and  allows decoupling when we deal with the Laplace transform of a sum of self adjoint random matrices. 
\begin{lemma} \label{lemme1} Assume that \eqref{condmajlamdamax} holds. Then for any  $t>0$ and any $k  \in \{0, \dots , \ell-1 \}$, 
\[
 \E {\rm Tr}\Big ({\rm e}^{t \sum_{j=1}^{2^k}  \bbS_{j}^{(k)}  } \Big )\le \E {\rm Tr}\Big ({\rm e}^{t \sum_{j=1}^{2^{k+1}}  \bbS_{j}^{(k+1)}   } \Big )  \Big ( 1+ \beta_{d_k + 1} {\rm e}^{tM n_{\ell} 2^{\ell - k }}\Big )^{2^k} \, ,
\]
where $(\bbS_{j}^{(k)} )_{j=1, \dots, 2^k}$ is the family of mutually independent  random matrices  defined in \eqref{1rstdef}. 
\end{lemma}

\noindent{\it Proof.} Note that for any $k  \in \{0, \dots , \ell-1 \}$ and any $j \in \{1, \dots, 2^k\}$, \[K_{k,j}= K_{k+1,2j-1} \cup K_{k+1,2j}  \] 
where the union is disjoint.  Therefore
\[
\bbS_{j}^{(k)}  = \bbS_{j,1}^{(k)}  + \bbS_{j,2}^{(k)}  \ \text{ and } \ \bV_j^{(k)} = ( \bV_{j,1}^{(k)} , \bV_{j,2}^{(k)} )  \, , 
\]
where $\bbS_{j,1}^{(k)}  := \sum_{i \in K_{k+1,2j-1}} \bbX_{j}^{(k)} (i)$, $\bbS_{j,2}^{(k)} := \sum_{i \in K_{k+1,2j}} \bbX_{j}^{(k)} (i)  $,   $\bV_{j,1}^{(k)} := \big (\bX_{j}^{(k)} (i)  \, , i \in K_{k+1, 2j-1} \big )$ and $\bV_{j,2}^{(k)} := \big (\bX_{j}^{(k)} (i)  \, , i \in K_{k+1, 2j} \big )$. 
Note that there are exactly $d_k$ integers between $K_{k+1,2j-1}$ and $K_{k+1,2j}  $ and that for any $i \in \{1, \ldots, 2^{k+1}\}$,
\[
\mathrm{Card}( K_{k+1,i}) = \mathrm{Card}( K_{k+1,1}) = 2^{\ell-(k+1)} n_{\ell} \, .
\]

Recall that the probability space is assumed to be large enough to contain a sequence $(\delta_i, \eta_i)_{i \in {\mathbb Z}}$ of iid random variables uniformly distributed on $[0,1]^2$   independent of the sequence $(\bX_i)_{ i \geqslant 0}$. Therefore according to the remark on the existence of the family $(\bV_j^{(m)})_{1 \leq  j \leq  2^m}$ made at the beginning of Section \ref{sectionfundamentallemma}, the sequence $(\eta_i)_{i \in {\mathbb Z}}$ is independent of $(\bV_j^{(m)})_{1 \leq  j \leq  2^m}$.  According to Lemma \ref{berbeelma}   there exists a random vector $\btV_{1,2}^{(k)}$ of size $d^2 \mathrm{Card}( K_{k+1,2})$ with the same law as $\bV_{1,2}^{(k)}$ that is measurable with respect to 
$\sigma(\eta_{1}) \vee \sigma( \bV_{1,1}^{(k)}) \vee \sigma( \bV_{1,2}^{(k)})$, independent of $\sigma( \bV_{1,1}^{(k)}) $ and such that
\begin{eqnarray} \label{coupling1}
{\mathbb P} (\btV_{1,2}^{(k)} \neq  \bV_{1,2}^{(k)}  )  = \beta \big ( \sigma( \bV_{1,1}^{(k)}), \sigma(\bV_{1,2}^{(k)}) \big )  \leq  \beta_{d_{k}+1} \nonumber \, ,
\end{eqnarray}
where the inequality comes from the fact that, by relation  \eqref{equivalentdefbeta}, the quantity $\beta \big ( \sigma( \bV_{1,1}^{(k)}), \sigma(\bV_{1,2}^{(k)}) \big )$ depends only on the joint distribution  of  $( \bV_{1,1}^{(k)} , \bV_{1,2}^{(k)})  $  and therefore, by \eqref{defofVjm}, 
\[
\beta \big ( \sigma( \bV_{1,1}^{(k)}), \sigma(\bV_{1,2}^{(k)}) \big )  = \beta \big ( \sigma( \bX_i \, , i \in K_{k+1,1}  ), \ \sigma( \bX_i \, , i \in K_{k+1,2}  ) \big )  \leq  \beta_{d_{k}+1}  \, .
\]
Note that by construction, $\btV_{1,2}^{(k)} $ is independent of $\sigma \big (  \bV_{1,1}^{(k)} ,  ( \bV_{j}^{(k)}  )_{j=2, \dots, 2^k}\big ) $.

\smallskip
\noindent
For any $i \in K_{k+1,2} $, let 
\[
\btX_{1,2}^{(k)} (i)= \pi_i^{(k+1)} (\btV_{1,2}^{(k)} ) \ \text{ and } \  \bbtS_{1,2}^{(k)} = \sum_{i \in K_{k+1,2}} \bbtX_{1,2}^{(k)} (i) \, ,
\]
where $\bbtX_{1,2}^{(k)} (i)$ is the $d \times d$ random  matrix associated with the random vector $\btX_{1,2}^{(k)} (i)$. 

\noindent
With the notations above, we have
\begin{multline} \label{l1p1}
\E {\rm Tr} \exp \Big ({t \sum_{j=1}^{2^k}  \bbS_{j}^{(k)}  } \Big ) 
  =
 \E \Big (\I_{\btV_{1,2}^{(k)} =  \bV_{1,2}^{(k)} } {\rm Tr} \exp \big({t \sum_{j=1}^{2^k}  \bbS_{j}^{(k)} }  \big)\Big )   
+ \E \Big (\I_{\btV_{1,2}^{(k)} \neq  \bV_{1,2}^{(k)} } {\rm Tr} \exp \big({t \sum_{j=1}^{2^k}  \bbS_{j}^{(k)} }  \big)\Big )  
 \\  
\leq   \E {\rm Tr}\exp \Big ({t \bbS_{1,1}^{(k)}  + t\bbtS_{1,2}^{(k)}  + t \!\sum_{j=2}^{2^k}  \bbS_{j}^{(k)}  }\Big )   
+\! \E \Big (\!\I_{\btV_{1,2}^{(k)} \neq  \bV_{1,2}^{(k)} } {\rm Tr} \exp \big({t\! \sum_{j=1}^{2^k}  \bbS_{j}^{(k)} }  \big)\!\Big ) \, . 
\end{multline}
(With usual convention, $\sum_{j=\ell}^{2^k}  \bbS_{j}^{(k)} $ is the null vector if $\ell > 2^k$). 
By Golden-Thompson inequality, we have
\[
{\rm Tr}\exp \Big ({t  \sum_{j=1}^{2^k}  \bbS_{j}^{(k)}   }\Big) \leq  {\rm Tr}\Big ( {\rm e}^{t  \bbS_{1}^{(k)} }\cdot {\rm e}^{t \sum_{j=2}^{2^k}  \bbS_{j}^{(k)}   }
  \Big ).
\]
Hence, since $\sigma \big (  \bV_{j}^{(k)}  \, , \, j=2, \dots, 2^k\big ) $ is independent of $\sigma  \big ( \bV_{1,1}^{(k)}  , \bV_{1,2}^{(k)} , \btV_{1,2}^{(k)}  \big )$, we get 
\begin{align*}
 \E \Big (\I_{\btV_{1,2}^{(k)} \neq  \bV_{1,2}^{(k)} } {\rm Tr} \exp\big({t \sum_{j=1}^{2^k}  \bbS_{j}^{(k)} } \big) \Big )
 \leq 
  {\rm Tr}\Big (\E \big( \I_{\btV_{1,2}^{(k)} \neq  \bV_{1,2}^{(k)} } {\rm e}^{t  \,  \bbS_{1}^{(k)}   }   \big)
   \cdot 
    \E  \exp\big({t \sum_{j=2}^{2^k}  \bbS_{j}^{(k)}   }\big)
  \Big )  \, .
\end{align*}
Note now the following fact: if ${\mathbb U} $ is a $d \times d$ self-adjoint random  matrix  with entries defined on $(\Omega, {\mathcal A}, {\mathbb P})$ and such that $\lambda_{{\rm max}} ({\mathbb U}) \leq b$ a.s., then for any $\Gamma \in {\mathcal A}$, 
\[
\I_\Gamma  {\rm e}^{{\mathbb U}} \preceq   {\rm e}^{b} {\mathbb I}_d \I_\Gamma  \  \text{ a.s. and so }  \ {\lambda}_{{\rm max}}  \E \big ( \I_\Gamma  {\rm e}^{{\mathbb U}}  \big ) \leq   {\rm e}^{b}  {\mathbb P} (\Gamma) \, .
\]
Therefore if we consider ${\mathbb V}$  a $d \times d$ self-adjoint random matrix  with entries defined on $(\Omega, {\mathcal A}, {\mathbb P})$, the following inequality is valid: 
\beq \label{tocomputeprodTrace}
{\rm Tr}  \big (  \E ( \I_\Gamma  {\rm e}^{{\mathbb U}} ) \E ( {\rm e}^{ {\mathbb V} } )\big ) \leq  {\rm e}^{b}  {\mathbb P} (\Gamma)  \cdot \E {\rm Tr}  (  {\rm e}^{ {\mathbb V} } )  \, .
\eeq
Notice now  that $ ( \bX_{1}^{(k)} (i) \, , i  \in K_{k,1})$ has the same distribution as $ ( \bX_i \, , i  \in K_{k,1})$. Therefore  $ \lambda_{{\rm max}} ( \bX_{1}^{(k)} (i)) \leq  M$ a.s. for any $i$, implying by Weyl's inequality that \[ \lambda_{{\rm max}} ( t \,   \bbS_{1}^{(k)}  ) \leq  t M  {\rm Card} (  K_{k,1} )  = t M 2^{\ell-k} n_{\ell} \ \text{ a.s.}\] Hence, applying \eqref{tocomputeprodTrace} with 
$ b = t M 2^{\ell-k} n_{\ell}$, $\Gamma= \{\btV_{1,2}^{(k)} \neq  \bV_{1,2}^{(k)} \}$ and ${\mathbb V} =  t  \sum_{j=2}^{2^k}  \bbS_{j}^{(k)}  $, and taking into account that ${\mathbb P} (\Gamma) \leq  \beta_{d_{k}+1}$, we obtain
\beq \label{beforeJensen}
  \E \Big (\I_{\btV_{1,2}^{(k)} \neq  \bV_{1,2}^{(k)} } {\rm Tr} \exp\big({t \sum_{j=1}^{2^k}  \bbS_{j}^{(k)} } \big) \Big )\!
  \leq 
   \beta_{d_{k}+1} {\rm e}^{  t n_{\ell} 2^{\ell  -k }  M }   \E {\rm Tr} \exp \Big ({t \sum_{j=2}^{2^k}  \bbS_{j}^{(k)} } \Big ) 
 \, . 
\eeq
Note now that if ${\mathbb V}$ and ${\mathbb W}$ are two independent random matrices with entries defined on $(\Omega, {\mathcal A}, {\mathbb P})$ and  such that $\E( {\mathbb W}) = {\mathbf 0}$ then 
\[
\E {\rm Tr}\exp ({\mathbb V}) = \E {\rm Tr}\exp \big ( \E \big (  {\mathbb V}  +  {\mathbb W} | \sigma( {\mathbb V}) \big ) \big )  \, .
\]
Since ${\rm Tr} \circ \exp $ is convex, it follows from Jensen's inequality applied to the conditional expectation that 
\beq \label{applicondJensen}
\E {\rm Tr}\exp ({\mathbb V})  \leq \E \big (  \E \big (  {\rm Tr} \,  {\rm e}^{  {\mathbb V } + {\mathbb W} } | \sigma( {\mathbb V}) \big )  \big ) =   \E \big (  {\rm Tr} \,  {\rm e}^{  {\mathbb V } + {\mathbb W} }  \big ) \, .
\eeq
Since  $\E ( \bbX_1^{(k)} (i) ) = \E ( \bbX_i) = {\mathbf 0}$ for any $i \in K_{k,1}$ and $\sigma(\bbS_{1,1}^{(k)} , \bbtS_{1,2}^{(k)} ) $ is independent of  $\sigma ( \bbS_j^{(k)}, j=2, \ldots ,2^k)$, we can apply the inequality above with ${\mathbb W}= t( \bbS_{1,1}^{(k)}  +\bbtS_{1,2}^{(k)} ) $ and ${\mathbb V} = t  \sum_{j=2}^{2^k}  \bbS_{j}^{(k)} $. Therefore, starting from \eqref{beforeJensen} and using \eqref{applicondJensen}, we get 
\beq \label{l1p2}
  \E \Big (\I_{\btV_{1,2}^{(k)} \neq  \bV_{1,2}^{(k)} } {\rm Tr} \exp\big({t \sum_{j=1}^{2^k}  \bbS_{j}^{(k)} } \big) \Big ) 
 \leq   \beta_{d_{k}+1} {\rm e}^{  t n_{\ell} 2^{\ell  -k }  M }  \, \E {\rm Tr} \exp\Big ({t \big(\bbS_{1,1}^{(k)}  + \bbtS_{1,2}^{(k)} +\sum_{j=2}^{2^k}  \bbS_{j}^{(k)} \big )} \Big)  \, .
 \eeq
Starting from \eqref{l1p1}  and considering \eqref{l1p2}, it follows that  
\beq \label{l1p3}
 \E {\rm Tr} \exp \Big ({t \sum_{j=1}^{2^k}  \bbS_{j}^{(k)}  } \Big ) \leq  (1+\beta_{d_{k}+1} {\rm e}^{  t n_{\ell} 2^{\ell  -k }  M })  \, \E {\rm Tr} \exp\Big ({t \big(\bbS_{1,1}^{(k)}  + \bbtS_{1,2}^{(k)} +\sum_{j=2}^{2^k}  \bbS_{j}^{(k)} \big )} \Big) \, .
 \eeq
The proof of Lemma \ref{lemme1} will then be achieved after having iterated this procedure $2^k-1$ times more. For the sake of clarity, let us explain the way to go from the $j$-th step to the $(j+1)$-th step. 

At the end of the $j$-th step,  assume  that we have constructed  with the help of the coupling Lemma  \ref{berbeelma}, $j$ random vectors  $\btV_{i,2}^{(k)}$, $i=1, \dots, j$,  each of dimension $d^2 \mathrm{Card}( K_{k+1,1})$ and satisfying the following properties:  for any $ i$  in $\{1, \dots, j \}$, $\btV_{i,2}^{(k)}$ is a measurable function of $( \bV_{i,1}^{(k)} , \bV_{i,2}^{(k)}, \eta_i )$, it has the same distribution as 
$\bV_{i,2}^{(k)}$, is such that ${\mathbb P} (\btV_{i,2}^{(k)} \neq  \bV_{i,2}^{(k)}  )  \leq  \beta_{d_{k}+1}$,  is    independent of $\bV_{i,1}^{(k)} $ and  it satisfies
\beq \label{l1p7}
\E {\rm Tr} \exp  \Big ( { t \sum_{j=1}^{2^k}  \bbS_{j}^{(k)}   }\Big ) 
\leq  \big ( 1+ \beta_{d_{k}+1} {\rm e}^{  t n_{\ell} 2^{\ell  -k }  M }  \big  )^j \cdot  \E {\rm Tr}  \exp \Big  (  {t  \sum_{i=1}^{j}  ( \bbS_{i,1}^{(k)}  +  \bbtS_{i,2}^{(k)}  ) + t\sum_{i=j+1}^{2^k}  \bbS_{i}^{(k)}   }\Big )  \, ,
\eeq 
where we have implemented the following notation: 
 \beq \label{notcoor}
\bbtS_{i,2}^{(k)} = \sum_{r \in K_{k+1,2i}} \bbtX_{i,2}^{(k)} (r) \, .
\eeq
In the notation above, $\bbtX_{i,2}^{(k)} (r) $ is the $d \times d$ random matrix associated with the random vector $\btX_{i,2}^{(k)} (r) $ of ${\mathbb K}^{d^2}$ defined by
\[
\btX_{i,2}^{(k)} (r)= \pi_r^{(k+1)} (\btV_{i,2}^{(k)} )   \text{ for any } r \in K_{k+1,2i}  \, .
\]
Note that the induction assumption above has been proven at the beginning of the proof  for $j=1$. Moreover, note that since, for any $m \in \{1, \ldots, \ell \} $, $(\bV_j^{(m)})_{1 \leq  j \leq  2^m}$ is a family of independent random vectors, the random vectors $\btV_{i,2}^{(k)}$, $i=1, \dots, j$, defined above are also such that, for any $ i \in \{1, \dots, j \}$,  $\btV_{i,2}^{(k)} $ is independent of $\sigma \big (  ( \bV_{\ell,1}^{(k)} )_{\ell=1, \dots, i}, 
 (\btV_{\ell,2}^{(k)} )_{\ell=1, \dots, i-1}  ,   ( \bV_{\ell}^{(k)}  )_{\ell=i+1, \dots, 2^k}\big ) $.

Now to show that the induction hypothesis also holds at step $j+1$, we proceed as follows. By  Lemma  \ref{berbeelma},  there exists a random vector $\btV_{j+1,2}^{(k)}$ of size $d^2 \mathrm{Card}( K_{k+1,1})$ with the same law as $\bV_{j+1,2}^{(k)}$,   measurable with respect to 
$\sigma(\eta_{j+1}) \vee \sigma( \bV_{j+1,1}^{(k)}) \vee \sigma( \bV_{j+1,2}^{(k)})$, independent of $\sigma( \bV_{j+1,1}^{(k)}) $ and such that
\begin{eqnarray} \label{coupling3}
{\mathbb P} (\btV_{j+1,2}^{(k)} \neq  \bV_{j+1,2}^{(k)}  )  \leq  \beta_{d_{k}+1}  \, .
\end{eqnarray}
(The inequality above comes again from \eqref{defofVjm} and the equivalent definition  \eqref{equivalentdefbeta} of the $\beta$-coefficients). Note that by construction, $\sigma \big (  ( \bV_{i,1}^{(k)} )_{i=1, \dots, j+1}, 
 (\btV_{i,2}^{(k)} )_{i=1, \dots, j}  ,   ( \bV_{i}^{(k)}  )_{i=j+2, \dots, 2^k}\big ) $    and  $\sigma(\btV_{j+1,2}^{(k)} )$ are independent. 
With the notation \eqref{notcoor}, we  have the following decomposition: 
\begin{multline}  \label{dec1stepjj+1}
 \E {\rm Tr}\exp \Big ( t  \sum_{i=1}^{j}  ( \bbS_{i,1}^{(k)}  +  \bbtS_{i,2}^{(k)} ) + t  \sum_{i=j+1}^{2^k}  \bbS_{i}^{(k)}   \Big )     
  \leq    
 \E {\rm Tr} \exp \Big ({t   \sum_{i=1}^{j+1}  ( \bbS_{i,1}^{(k)}  +  \bbtS_{i,2}^{(k)} ) + t \sum_{i=j+2}^{2^k}  \bbS_{i}^{(k)}   }\Big )  
 \\   + 
 \E \Big( \I_{\btV_{j+1,2}^{(k)} \neq  \bV_{j+1,2}^{(k)} }  {\rm Tr}\exp \big ( {t  \sum_{i=1}^{j}  ( \bbS_{i,1}^{(k)}  + \bbtS_{i,2}^{(k)} ) +t \sum_{i=j+1}^{2^k}  \bbS_{i}^{(k)}   }  \big)\Big )  \, .
\end{multline}
Using Golden-Thompson inequality, we have
\begin{align*}
{\rm Tr}\exp \Big (& {t  \sum_{i=1}^{j}  ( \bbS_{i,1}^{(k)}  +  \bbtS_{i,2}^{(k)} ) +t \sum_{i=j+1}^{2^k}  \bbS_{i}^{(k)}   }  \Big)
\\ &
 \leq  {\rm Tr} \Big (\exp \big( t \, \bbS_{j+1}^{(k)}  \big) \cdot \exp \big( t \sum_{i=1}^{j}  ( \bbS_{i,1}^{(k)}  + \bbtS_{i,2}^{(k)} ) + t\sum_{i=j+2}^{2^k}  \bbS_{i}^{(k)} \big) \Big )  \, .
\end{align*}
Hence, since the sigma algebra generated by  $ \big (  ( \bV_{i,1}^{(k)} )_{i=1, \dots, j}, 
 (\btV_{i,2}^{(k)} )_{i=1, \dots, j}  ,   ( \bV_{i}^{(k)}  )_{i=j+2, \dots, 2^k}\big ) $ is independent of that generated by $  \big ( \bV_{j+1,1}^{(k)}  , \bV_{j+1,2}^{(k)} , \btV_{j+1,2}^{(k)}  \big )$, we get 
\begin{multline} \label{stepjtrace}
 \E {\rm Tr}\Big (  {\rm e}^{t \big ( \sum_{i=1}^{j}  ( \bS_{i,1}^{(k)}   +  \btS_{i,2}^{(k)} ) + \sum_{i=j+1}^{2^k}  \bS_{i}^{(k)}  \big ) }  \I_{\btV_{j+1,2}^{(k)} \neq  \bV_{j+1,2}^{(k)} } \Big ) \\
  \leq  {\rm Tr} \Big ( \E {\rm e}^{t \big ( \sum_{i=1}^{j}  ( \bS_{i,1}^{(k)}   +  \btS_{i,2}^{(k)} ) + \sum_{i=j+2}^{2^k}  \bS_{i}^{(k)}  \big ) }   
 \cdot \E \big (  {\rm e}^{t  \,  \bS_{j+1}^{(k)}   }  \I_{\btV_{j+1,2}^{(k)} \neq  \bV_{j+1,2}^{(k)} }  \big ) \Big ) \, .
 \end{multline}
 By Weyl's inequality,
 \[
 \lambda_{{\rm max}}  \big ( t \,   \bS_{j+1}^{(k)}    \big ) \leq  t \sum_{r \in K_{k,j+1}}   \lambda_{{\rm max}}   ( \bbX_{j+1}^{(k)} (r) ) \ \text{a.s.} 
 \]
Using that  $\bV_{j+1}^{(k)} =^{\mathcal D} (\bX_i \, , i \in K_{k,j+1} )$ and that $ \lambda_{{\rm max}} (\bbX_i) \leq M$ a.s. for any $i$, it follows that 
\[
\lambda_{{\rm max}}  \big ( t \,   \bS_{j+1}^{(k)}      \big ) \leq    tM  {\rm Card} (K_{k,j+1} )  = tM2^{\ell -k} n_k  \ \text{a.s.} 
\]
In addition, we notice that $ \bS_{j+1,1}^{(k)}   +  \btS_{j+1,2}^{(k)}  $ is independent of $\sum_{i=1}^{j}  ( \bS_{i,1}^{(k)}   +  \btS_{i,2}^{(k)} ) + \sum_{i=j+2}^{2^k}  \bS_{i}^{(k)} $ and, since $\btV_{j+1,2}^{(k)}=^{{\mathcal D}}\bV_{j+1,2}^{(k)}$ and  $\bV_{j+1}^{(k)} =^{\mathcal D} (\bX_i \, , i \in K_{k,j+1} )$, $\E( \bS_{j+1,1}^{(k)}   +  \btS_{j+1,2}^{(k)} ) ={\bf 0}$. 
Therefore, starting from \eqref{stepjtrace} and taking into account \eqref{coupling3}, an application of  inequality \eqref{tocomputeprodTrace} with $ b = t M 2^{\ell-k} n_{\ell}$, $\Gamma= \{\btV_{j+1,2}^{(k)} \neq  \bV_{j+1,2}^{(k)} \}$ and ${\mathbb V} =  t \big ( \sum_{i=1}^{j}  ( \bS_{i,1}^{(k)}   +  \btS_{i,2}^{(k)} ) + \sum_{i=j+2}^{2^k}  \bS_{i}^{(k)}  \big ) $, followed by an application of inequality \eqref{applicondJensen} with ${\mathbb W} = t ( \bS_{j+1,1}^{(k)}   +  \btS_{j+1,2}^{(k)} )$,  gives
\begin{multline} \label{dec2stepjj+1}
 \E {\rm Tr}\Big (  {\rm e}^{t \big ( \sum_{i=1}^{j}  ( \bS_{i,1}^{(k)}   +  \btS_{i,2}^{(k)} ) + \sum_{i=j+1}^{2^k}  \bS_{i}^{(k)}  \big ) }  \I_{\btV_{j+1,2}^{(k)} \neq  \bV_{j+1,2}^{(k)} } \Big )
 \\ \leq  \beta_{d_k+1} {\rm e}^{  t n_{\ell} 2^{\ell  -k }  M }   \times \E {\rm Tr}\Big ({\rm e}^{t  \big ( \sum_{i=1}^{j+1}  ( \bS_{i,1}^{(k)}  +  \btS_{i,2}^{(k)} ) + \sum_{i=j+2}^{2^k}  \bS_{i}^{(k)}  \big ) }\Big )  \, .
 \end{multline}
Therefore, starting from \eqref{dec1stepjj+1} and using \eqref{dec2stepjj+1}, we get 
\begin{align*}
\E {\rm Tr} \Big  ( &  {\rm e}^{t \big ( \sum_{i=1}^{j}  ( \bS_{i,1}^{(k)}  +  \btS_{i,2}^{(k)}  ) + \sum_{i=j+1}^{2^k}  \bS_{i}^{(k)}  \big ) }\Big )  \nonumber \\
& \leq
 \big ( 1+ \beta_{d_k+1} {\rm e}^{  t n_{\ell} 2^{\ell  -k }  M }  \big  ) \times \E {\rm Tr}\Big ({\rm e}^{t  \big ( \sum_{i=1}^{j+1}  ( \bS_{i,1}^{(k)}  +  \btS_{i,2}^{(k)} ) + \sum_{i=j+2}^{2^k}  \bS_{i}^{(k)}  \big ) }\Big )    \, ,
\end{align*}
which combined with \eqref{l1p7} implies that 
\[
\E {\rm Tr} \Big ({\rm e}^{t   \sum_{j=1}^{2^k}  \bS_{j}^{(k)}   } \Big )  \leq    \big ( 1+ \beta_{d_k+1} {\rm e}^{  t n_{\ell} 2^{\ell  -k }  M }  \big  )^{j+1} \times  \E {\rm Tr} \Big ({\rm e}^{t  \big ( \sum_{i=1}^{j+1}  ( \bS_{i,1}^{(k)}  +  \btS_{i,2}^{(k)} ) + \sum_{i=j+2}^{2^k}  \bS_{i}^{(k)}  \big ) }\Big )   \, ,
\]
proving the induction hypothesis for the step $j+1$.  Finally $2^k$ steps of the procedure lead to 
\beq \label{l1p8}
\E {\rm Tr}\Big ({\rm e}^{t   \sum_{j=1}^{2^k}  \bS_{j}^{(k)}  } \Big )  \leq    \big ( 1+ \beta_{d_k+1} {\rm e}^{  t n_{\ell} 2^{\ell  -k }  M }  \big  )^{2^k} \times  \E {\rm Tr} \Big ({\rm e}^{t    \sum_{i=1}^{2^k}  ( \bS_{i,1}^{(k)}  +  \btS_{i,2}^{(k)} )  }\Big )   \, .
\eeq
To end the proof of the lemma it suffices to notice the following facts: the random vectors $ \bV_{i,1}^{(k)}  ,  \btV_{i,2}^{(k)} $, $i=1, \ldots, 2^k$, are mutually independent and such that 
$ \bV_{i,1}^{(k)} =^{\mathcal D}  \bV_{2i-1}^{(k+1)}$ and $ \bV_{i,2}^{(k)} =^{\mathcal D}  \bV_{2i}^{(k+1)}$. In addition, the random vectors $ \bV_{i}^{(k+1)}$, $i=1, \ldots, 2^{k+1}$, are mutually independent. This obviously implies that 
\[
\E {\rm Tr} \Big ({\rm e}^{t    \sum_{i=1}^{2^k}  ( \bS_{i,1}^{(k)}  +  \btS_{i,2}^{(k)} )  }\Big ) = \E {\rm Tr} \Big ({\rm e}^{t    \sum_{i=1}^{2^{k+1}}   \bS_{i}^{(k+1)}  }\Big ) \, ,
\]
which ends the proof of the lemma.  $\square$

\medskip

\subsubsection{Proof of Proposition \ref{propinter1}}
We shall prove Inequality \eqref{resultpropinter1} with $K_{A}$ defined in Section \ref{sectionCantor}.

\smallskip

Let us  prove it first  in the case where  $0 < tM \leq  4/A$. Since by  Weyl's inequality, 
\[
 \lambda_{{\rm max}} \big ( \sum_{i \in K_{A}}\bbX_i  \big )  \leq  \sum_{i \in K_{A}}  \lambda_{{\rm max}} \big (\bbX_i  \big )  \leq  A M \, ,
\]
and $\E ( {\mathbb X}_i) = {\mathbf 0}$ for any $i \in K_A$, it follows by using Lemma \ref{lmatropp} applied with $K=\{1\}$ and ${\mathbb U}_1=\sum_{i \in K_{A}}\bbX_i$ that, for any $t >0$,
\[
\E {\rm Tr}\big ({\rm e}^{t  \sum_{i \in K_{A}}\bbX_i   } \big )   \leq  d  \exp  \Big (    t^2 g(t AM) \lambda_{{\rm max}}  \Big (\E \big (  \sum_{i \in K_A} {\mathbb X}_i \big )^2 \Big )  \Big )\, .
\]
Therefore by the definition of $v^2$, since $g$ is increasing, $tAM < 4$ and $g(4) \leq  3.1$, we get 
\[
\E {\rm Tr}\big ({\rm e}^{t  \sum_{i \in K_{A}}\bbX_i   } \big )   \leq  d  \exp   (  3.1 \times  A t^2 v^2 ) \, ,  
\]
proving then \eqref{resultpropinter1}. 

\smallskip

We prove now Inequality \eqref{resultpropinter1} in the case where  $4/A < tM \leq \min  \big (  \frac{1}{2}, \frac{c \log 2}{32 \log A} \big )$. Let 
\begin{equation}
\kappa = \frac{c}{8}  \ \text{ and } \ k(t)=\inf \Big \{k\in \mathbb{Z}:A((1-\delta )/2)^{k}\leq 
\min \big ( \frac{\kappa}{(tM)^{2}} , A \big ) \Big \}. \label{defk_2}
\end{equation}
Note that if $t^{2}M^{2}\leq  \kappa /A$ then $k(t) =0$ whereas $k(t) \geq 1$ if  $t^{2}M^{2}>\kappa  /A$. In addition  by the selection of $\ell_{A}$,  $A((1-\delta )/2)^{\ell} < 4/\delta$. Therefore  $k  (t)\leq  \ell_{A}$ since 
$(tM)^2\leq c \delta/32 $. Then, starting from \eqref{trivial}, considering the selection of $k(t)$ and  using  Lemma \ref{lemme1}, we get by induction that 
\beq \label{constrivial2} 
\E {\rm Tr}\exp \big (  t   \sum_{i \in K_{A}} \bbX_i    \big )  
\leq  \prod_{k=0}^{k(t) -1} \Big ( 1+ \beta_{d_{k}+1 }{\rm e}^{tM n_{\ell} 2^{\ell - k }}\Big )^{2^k}\, \E {\rm Tr} \exp \Big (t \sum_{j=1}^{2^{ k(t)}}  \bS_{j}^{(k(t))} \Big )  \, ,
\eeq 
with the usual convention that $\prod_{k=0}^{-1}  a_k = 1$. Note that in the inequality 
above,  $(\bbS_{j}^{(k(t))} )_{j=1, \ldots, 2^{ k(t)}}$ is a family of mutually independent  random matrices  defined in \eqref{1rstdef}. They are then constructed from a family $(\bV_j^{(k(t))})_{1 \leq  j \leq  2^{k(t)}}$  of  
$2^{k(t)}$ mutually independent random vectors that satisfy \eqref{defofVjm}. Therefore we have that, for any $ j \in 1, \ldots, 2^{ k(t)}$, $\bbS_{j}^{(k(t))} =^{\mathcal D} \sum_{i \in K_{k(t),j}} \bbX_i $. Moreover, according to the remark on the existence of the family $(\bV_j^{(k(t))})_{1 \leq  j \leq  2^{k(t)}}$ made at the beginning of Section \ref{sectionfundamentallemma},  the entries of each random matrix $\bbS_{j}^{(k(t))} $ are measurable functions of $(\bX_i , \delta_i)_{i \in {\mathbb Z}}$. 

The rest of the proof consists of giving  a suitable upper bound for $\E {\rm Tr} \exp \Big (t \sum_{j=1}^{2^{ k(t)}}  \bS_{j}^{(k(t))} \Big ) $. With this aim,  let $p$ be a positive integer to be chosen later  such that 
\beq \label{contraintep}
2p \leq  {\rm Card}   (  K_{k(t),j} ) :=q \, .
\eeq
Note that  $q= 2^{\ell -k(t) } n_\ell $ and by \eqref{evident} 
\[
q \geq \frac{A}{2^{k(t) +1} } \, .
\]
Therefore if $k(t) =0$ then $q \geq A/2$ implying that $q \geq 2$ (since we have $4/A < tM \leq 1$). Now if $k(t) \geq 1$ and therefore if $t^{2}M^{2}>\kappa /A$, by the definition of $k(t)$, we have $q \geq \frac{\kappa}{(tM)^{2}}$  and then $q \geq 2$ since $ (tM)^2 \leq \kappa/2$.  Hence in all cases, $q \geq 2$ implying that the selection of a positive integer $p$ satisfying \eqref{contraintep} is always possible.

Let $m_{q,p} = [q/(2p)]$. For any $j  \in \{1, \dots , 2^{k(t)} \}$, we divide $K_{k(t), j} $ into 
$2m_{q,p} $ consecutive
intervals $(J^{( k(t) )}_{j,i} \, , \, 1 \leq  i \leq  2m_{q,p}  )$ each containing  $p$ consecutive integers plus a remainder interval $J^{( k(t) )}_{j, 2 m_{q,p}  +1}  $ containing $r = q  - 2 p m_{q,p}  $ consecutive integers. Note that 
this last interval contains at most $2p-1$ integers.  Let  $\bbX_{j}^{(k(t))} (k) $ be the $d \times d$ random matrix  associated with the random vector   $\bX_{j}^{(k(t))} (k) $ defined in \eqref{1rstdef} and define
\beq \label{defbZ}
\bbZ_{j,i  }^{(k(t))} = \sum_{k \in K_{k(t),j} \cap J^{( k(t) )}_{j,i}} \bbX_{j}^{(k(t))} (k)  \, .
\eeq
With this notation
\[
\bS_{j  }^{(k(t))}  = \sum_{ i=1}^{m_{q,p} +1} \bbZ_{j , 2i -1}^{(k(t))}  + \sum_{ i=1}^{m_{q,p} } \bbZ_{j , 2i }^{(k(t))}   \, .
\]
Since ${\rm Tr} \circ \exp$ is a convex function, we get 
\begin{equation}\label{trace-exp-convex}
\E {\rm Tr} \exp \Big ( t\sum_{j=1}^{2^{ k(t)}}  \bS_{j}^{(k(t))} \Big ) 
\leq  \frac{1}{2} \E {\rm Tr} \exp \Big (2t  \sum_{j=1}^{2^{ k(t) }}  \sum_{ i=1}^{m_{q,p} +1} \bbZ_{j , 2i -1}^{(k(t))}  \Big )
+ \frac{1}{2} \E {\rm Tr} \exp \Big (2 t\sum_{j=1}^{2^{ k(t) }}  \sum_{ i=1}^{m_{q,p} } \bbZ_{j , 2i }^{(k(t))}  \Big ) \, .
\end{equation}
We start by giving an upper bound for $ \E {\rm Tr} \exp \Big (2 t  \sum_{j=1}^{2^{ k(t) }}  \sum_{ i=1}^{m_{q,p} } \bbZ_{j , 2i }^{(k(t))}   \Big ) $. With this aim, let us define the following vectors
\beq \label{defbU}
\bU_{j ,i }^{(k(t))}  = \big ( \bX_{j}^{(k(t))} (k) \, , \, k \in K_{k(t),j} \cap J^{( k(t) )}_{j,i}  \big  )  \, \text{ and } \,  {\mathbf W}_j^{(k(t))}  = \big  ( \bU_{j ,i }^{(k(t))}  \, , \, i \in \{1, \ldots, 2m_{q,p} +1\} \big )    \, .
\eeq
Proceeding by induction and using  the coupling lemma \ref{berbeelma}, one can construct  random vectors  $\bU_{j,2i }^{*(k(t))}  $, $j=1, \dots, 2^{k(t)} $, $i=1, \dots, m_{q,p}$, that satisfy the following properties:
\begin{enumerate}[label=(\roman{*})]
\item $(\bU_{j,2i }^{*(k(t))}, (j,i) \in \{ 1, \dots, 2^{k(t)}\} \times \{ 1, \dots, m_{q,p} \} )$ is a family of  mutually independent random vectors,
\item $ \bU_{j,2i }^{*(k(t))}  $ has the same distribution as $ \bU_{j,2i }^{(k(t))}  $,
\item \label{bloccoupling}$ {\mathbb P} (\bU_{j ,2i}^{*(k(t))}   \neq \bU_{j,2i }^{(k(t))}   )  \leq  \beta_{p+1}  \, . $
\end{enumerate}
Let us explain the construction. Recall first that $(\Omega,{\cal A}, {\mathbb P} )$  is assumed to be rich enough  to contain a sequence  $ (\eta_i)_{i \in {\mathbb Z}}$ of iid random variables with uniform distribution over $[0,1]$ independent of $ (\bX_i, \delta_i)_{i \in {\mathbb Z}}$ (the sequence $(\delta_i)_{i \in {\mathbb Z}}$ has been used to construct the independent random matrices  $ \bS_{j}^{(k(t))}$, $j=1, \dots, k(t)$, involved in inequality \eqref{constrivial2}).  For any $j \in  \{ 1, \dots, 2^{k(t)}\} $, let  $\bU_{j ,2}^{*(k(t))}  = \bU_{j,2 }^{(k(t))}$, and construct the  random vectors $\bU_{j,2i }^{*(k(t))} $, $i =2, \dots,  m_{q,p}$, recursively from  $ ( \bU_{j,2\ell }^{*(k(t))} , 1 \leq \ell \leq i-1) $ as follows. According to Lemma \ref{berbeelma}, there exists a random vector $\bU_{j , 2i}^{*(k(t))} $  such that 
\beq \label{constrctionU2i}
\bU_{j ,2i}^{*(k(t))} =f_{i , j } \big (  (  \bU_{j,2\ell }^{*(k(t))})_{1 \leq \ell \leq i-1} , 
\bU_{j,2 i  }^{(k(t))} , \eta_{i + (j-1)2^{k(t)}} \big ) 
\eeq
 where $f_{i,j}$ is a measurable function, $\bU_{j ,2i}^{*(k(t))}$ has the same law as $\bU_{j ,2i}^{(k(t))} $, is independent of 
 $\sigma \big (    \bU_{j,2\ell }^{*(k(t))}  , 1 \leq \ell \leq i-1 \big ) $ and 
\[
 {\mathbb P} (\bU_{j,2i}^{*(k(t))}   \neq \bU_{j,2i }^{(k(t))}   )  = \beta \big ( \sigma \big (    \bU_{j,2\ell }^{*(k(t))}  , 1 \leq \ell \leq i-1 \big )   , \sigma( \bU_{j,2i }^{(k(t))}) \big ) \leq  \beta_{p+1} \, .
\]
By construction, for any fixed $j \in  \{ 1, \dots, 2^{k(t)}\} $, the random vectors $\bU_{j,2i }^{*(k(t))} $, $i =1, \dots,  m_{q,p}$, are mutually independent. In addition, by \eqref{constrctionU2i} and the fact that  $ ( {\mathbf W}_j^{(k(t))} , j = 1, \ldots, 2^{k(t)}  ) $ is a family of mutually independent  random vectors, we note that $(\bU_{j,2i }^{*(k(t))} , (i,j)\in \{1, \dots,  m_{q,p}\} \times \{1, \ldots, 2^{k(t)}  \}) $ is also so. Therefore the constructed random vectors $\bU_{j,2i }^{*(k(t))} $ $i =1, \dots,  m_{q,p}$, $j=1, \ldots, 2^{k(t)} $, satisfy Items {\rm (i)} and  {\rm (ii)} above. Moreover, by \eqref{constrctionU2i}, we have 
\[
\sigma \big (    \bU_{j,2\ell }^{*(k(t))}  , 1 \leq \ell \leq i-1 \big ) \subseteq \sigma \big (    \bU_{j,2\ell }^{(k(t))}   , 1 \leq \ell \leq i-1 \big ) \vee  \sigma \big (    \eta_{\ell + (j-1)2^{k(t)}}    , 1 \leq \ell \leq i-1 \big ) \, .
\]
Since $ (\eta_i)_{i \in {\mathbb Z}}$ is  independent of $ (\bX_i, \delta_i)_{i \in {\mathbb Z}}$, we have 
\[
 \beta \big ( \sigma \big (    \bU_{j,2\ell }^{*(k(t))}  , 1 \leq \ell \leq i-1 \big )   , \sigma( \bU_{j,2i }^{(k(t))}) \big )  \leq  \beta \big ( \sigma \big (    \bU_{j,2\ell }^{(k(t))}  , 1 \leq \ell \leq i-1 \big )   , \sigma( \bU_{j,2i }^{(k(t))}) \big )  \, .
\]
By relation  \eqref{equivalentdefbeta}, the quantity $\beta \big ( \sigma \big (    \bU_{j,2\ell }^{(k(t))}  , 1 \leq \ell \leq i-1 \big )   , \sigma( \bU_{j,2i }^{(k(t))}) \big )$ depends only on the joint distribution  of  $ \big (   (  \bU_{j,2\ell }^{(k(t))} )_{  1 \leq \ell \leq i-1 }   ,  \bU_{j,2i }^{(k(t))} \big )  $.  By the definition \eqref{defbU} of the $\bU_{j, \ell }^{(k(t))}$'s, the definition \eqref{1rstdef} of the $\bX_j^{(k(t))}(k) $'s and \eqref{defofVjm}, we infer that 
\begin{multline*}
\beta \big ( \sigma \big (    \bU_{j,2\ell }^{(k(t))}  , 1 \leq \ell \leq i-1 \big )   , \sigma( \bU_{j,2i }^{(k(t))}) \big ) \\
= \beta \big ( \sigma \big (   \bX_{k}\, , \, k \in \cup_{\ell=1}^{i-1}K_{k(t),j} \cap J^{( k(t) )}_{j,2 \ell} \big )   , \sigma \big (   \bX_{k}\, , \, k \in K_{k(t),j} \cap J^{( k(t) )}_{j,2 i} \big )  \big )
\leq \beta_{p+1}\, .
\end{multline*}
So, overall, the constructed random vectors $\bU_{j,2i }^{*(k(t))} $ $i =1, \dots,  m_{q,p}$, $j=1, \ldots, 2^{k(t)} $, satisfy also Item {\rm (iii)}  above.

Denote now 
\[
  \bX_{j,2i }^{*(k(t))}   ( \ell) = \pi_{\ell} (   \bU_{j,2i}^{*(k(t))} )
\]
where $\pi_i^{(m)}$ is the $\ell$-th canonical projection from ${\mathbb K}^{p d^2}$ onto ${\mathbb K}^{d^2}$, namely: for any vector ${\mathbf x} = ({\mathbf x}_i \, , i \in \{1, \ldots , p \}) $ of ${\mathbb K}^{p d^2}$,  $\pi_{\ell}( {\mathbf x} ) = {\mathbf x}_{\ell}$. Let  $\bbX_{j,2i }^{*(k(t))}   ( \ell) $ be the $d \times d$ random  matrix  associated with $\bX_{j,2i }^{*(k(t))}   ( \ell) $ and define, for any $i = 1 \ldots, m_{q,p}$, 
\[
{\mathbb Z}_{j,2i }^{*(k(t))}   = \sum_{\ell \in  K_{k(t),j} \cap J^{( k(t) )}_{j,2 i}}\bbX_{j,2i }^{*(k(t))}   ( \ell) \, .
\]
Observe that by Item {\rm (ii)} above, ${\mathbb Z}_{j,2i }^{*(k(t))}  =^{\mathcal D} {\mathbb Z}_{j,2i }^{(k(t))}  $ (where we recall that ${\mathbb Z}_{j,2i }^{(k(t))} $ is defined by \eqref{defbZ}) and by  Item {\rm (i)}, the random matrices ${\mathbb Z}_{j,2i }^{*(k(t))}   $,   $i =1, \dots,  m_{q,p}$, $j=1, \ldots, 2^{k(t)} $, are mutually independent. The aim now is to prove that the following inequality is valid: 
\beq\label{p1prop1}
\E   {\rm Tr}\exp \Big ( { 2t    \sum_{j=1}^{2^{ k(t) }}  \sum_{ i=1}^{m_{q,p} } {\mathbb Z}_{j,2i }^{(k(t))}    } \Big ) 
  \leq   \Big ( 1+( m_{q,p} -1)  {\rm e}^{qtM} \beta_{p+1} \Big )^{2^{k(t)}  } 
  \E {\rm Tr}  \exp \Big (2t \sum_{j=1}^{2^{k(t)}  } \sum_{ i=1}^{m_{q,p} } {\mathbb Z}_{j,2i }^{*(k(t))}   \Big ) \, . 
 \eeq
 Obviously, this can be done by induction if we can show that, for any $\ell$ in $\{1, \ldots, 2^{k(t)} \}$,
 \begin{multline}  \label{p1prop1induction}
\E   {\rm Tr}\exp \Big ( { 2t    \sum_{j=1}^{ \ell -1 }  \sum_{ i=1}^{m_{q,p} } {\mathbb Z}_{j,2i }^{*(k(t))}    + 2t    \sum_{j=\ell}^{2^{ k(t) }}  \sum_{ i=1}^{m_{q,p} }{\mathbb Z}_{j,2i }^{(k(t))}    } \Big ) \\
  \leq   \Big ( 1+ ( m_{q,p} -1)  {\rm e}^{qtM} \  \beta_{p+1} \Big )
\E   {\rm Tr}\exp \Big ( { 2t    \sum_{j=1}^{ \ell  }  \sum_{ i=1}^{m_{q,p} } {\mathbb Z}_{j,2i }^{*(k(t))}    + 2t    \sum_{j=\ell +1}^{2^{ k(t) }}  \sum_{ i=1}^{m_{q,p} }{\mathbb Z}_{j,2i }^{(k(t))}    } \Big )  \, .  
\end{multline}
To prove the inequality above, we set \[{\mathbb C}_{\ell-1, \ell}(t) = 2t    \sum_{j=1}^{ \ell -1 }  \sum_{ i=1}^{m_{q,p} } {\mathbb Z}_{j,2i }^{*(k(t))}    + 2t    \sum_{j=\ell}^{2^{ k(t) }}  \sum_{ i=1}^{m_{q,p} }{\mathbb Z}_{j,2i }^{(k(t))}    \] 
and we write
\begin{multline} \label{p1bisprop1induction}
\E   {\rm Tr}\exp \big ( {\mathbb C}_{\ell-1, \ell}(t) \big )  
  =    \E \Big (\prod_{i=2}^{m_{q,p}}\I_{\bU_{\ell,2i}^{(k(t))} =  \bU_{\ell,2i}^{*(k(t))} }  {\rm Tr}\exp \big ( {\mathbb C}_{\ell-1, \ell}(t) \big )  \Big ) 
 \\    \quad \quad +  \E \Big (\I_{ \exists i \in \{2, \dots, m_{q,p} \} \, : \, \bU_{\ell,2i}^{(k(t))} \neq  \bU_{\ell,2i}^{*(k(t))} }  {\rm Tr}\exp \big ( {\mathbb C}_{\ell-1, \ell}(t) \big )   \Big )  \\ 
 \leq    \E   {\rm Tr}\exp \big ( {\mathbb C}_{\ell, \ell +1}(t) \big )  +   \E \Big (\I_{ \exists i \in \{2, \dots, m_{q,p} \} \, : \, \bU_{\ell,2i}^{(k(t))} \neq  \bU_{\ell,2i}^{*(k(t))} }  {\rm Tr}\exp \big ( {\mathbb C}_{\ell-1, \ell}(t) \big )   \Big )    \, .
\end{multline}
Note that the sigma algebra generated by the random vectors $( \bU_{j,2i}^{*(k)}  )_{i \in \{1, \ldots, m_{q,p} \} , j \in \{1, \ldots, \ell-1\}}$ and $( \bU_{j,2i}^{(k)}  )_{i \in \{1, \ldots, m_{q,p} \}, \,  j \in \{\ell +1, \ldots, 2^{ k(t) }\}}$ is independent of $ \sigma \big ( (  \bU_{\ell,2i}^{(k)}, \bU_{\ell,2i}^{*(k)}  )_{i \in \{1, \dots, m_{q,p} \}}  \big )$. This fact  together with the 
Golden-Thomson inequality give
\begin{align*}
 \E &  \Big (\I_{ \exists i \in \{2, \dots, m_{q,p} \} \, : \, \bU_{\ell,2i}^{(k(t))} \neq  \bU_{\ell,2i}^{*(k(t))} }  {\rm Tr}\exp \big ( {\mathbb C}_{\ell-1, \ell}(t) \big )   \Big ) 
\\
&\leq 
{\rm Tr} \Big ( \E \Big( \exp \big (2t  \sum_{j=1}^{ \ell -1 }  \sum_{ i=1}^{m_{q,p} } \bbZ_{j , 2i }^{*(k(t))} + 2t \sum_{j=\ell +1}^{2^{ k(t) }}  \sum_{ i=1}^{m_{q,p} } \bbZ_{j , 2i }^{(k(t))}  \big )\Big)  
\\
 & \quad  \quad \quad \quad  \quad \quad \quad  \quad \quad \times \,  
\E \Big (\I_{ \exists i \in \{2, \dots, m_{q,p} \} \, : \, \bU_{\ell,2i}^{(k(t))} \neq  \bU_{\ell,2i}^{*(k(t))} }    \exp \big(2t  \sum_{ i=1}^{m_{q,p} } \bbZ_{\ell, 2i }^{(k(t))} \big) \Big)   \Big )  \, .
 \end{align*}
By  Weyl's  inequality and \eqref{defofVjm}, we infer that, almost surely, 
\beq \label{lambdamaxpair}
\lambda_{{\rm max}}  \big(2t  \sum_{ i=1}^{m_{q,p} } \bbZ_{\ell, 2i }^{(k(t))}   \big) 
\leq  2t \sum_{ i=1}^{m_{q,p} }  \sum_{k \in K_{k(t),\ell} \cap J^{( k(t) )}_{\ell ,2i}}  \lambda_{{\rm max}}\big (  \bbX_k  \big ) 
\leq  2tm_{q,p}  p M 
\leq  tqM \, .
\eeq
Therefore, applying \eqref{tocomputeprodTrace} with $ b =   tqM $, $\Gamma= \{ \exists i \in \{2, \dots, m_{q,p} \} \, : \, \bU_{\ell,2i}^{(k(t))} \neq  \bU_{\ell,2i}^{*(k(t))} \}$ and ${\mathbb V} = 2t  \sum_{j=1}^{ \ell -1 }  \sum_{ i=1}^{m_{q,p} } \bbZ_{j , 2i }^{*(k(t))} + 2t \sum_{j=\ell +1}^{2^{ k(t) }  }  \sum_{ i=1}^{m_{q,p} } \bbZ_{j , 2i }^{(k(t))} $ and taking into account that 
\[
{\mathbb P}  (    \Gamma ) \leq  \sum_{i=2}^{ m_{q,p}}{\mathbb P} ( \bU_{\ell,2i}^{(k(t))} \neq  \bU_{\ell,2i}^{*(k(t))}   )  \leq    ( m_{q,p}  -1)  \beta_{p+1}  \, ,
\]
we get 
\begin{multline*}
 \E   \Big (\I_{ \exists i \in \{2, \dots, m_{q,p} \} \, : \, \bU_{\ell,2i}^{(k(t))} \neq  \bU_{\ell,2i}^{*(k(t))} }  {\rm Tr}\exp \big ( {\mathbb C}_{\ell-1, \ell}(t) \big )   \Big )    \\
\leq 
( m_{q,p} -1)     \beta_{p+1} {\rm e}^{qtM}    \, \E {\rm Tr}  \exp \Big ( 2t  \sum_{j=1}^{ \ell -1 }  \sum_{ i=1}^{m_{q,p} } \bbZ_{j , 2i }^{*(k(t))} + 2t \sum_{j=\ell +1}^{2^{ k(t) }  }  \sum_{ i=1}^{m_{q,p} } \bbZ_{j , 2i }^{(k(t))}  \Big ) \, .
\end{multline*}
Using   that the sigma algebra generated by the random vectors $( \bU_{j,2i}^{*(k)}  )_{i \in \{1, \ldots, m_{q,p} \} , \,  j \in \{1, \ldots, \ell-1\}}$ and $( \bU_{j,2i}^{(k)}  )_{i \in \{1, \ldots, m_{q,p} \}, j \in \{\ell +1, \ldots, 2^{ k(t) }\}}$ is independent of $ \sigma \big ( ( \bU_{\ell,2i}^{*(k)}  )_{i \in \{1, \dots, m_{q,p} \}}  \big )$, and noticing that by construction, 
$\E (\bbZ_{\ell , 2i }^{*(k(t))}) = \E (\bbZ_{\ell , 2i }^{(k(t))})={\bf 0}$, an application of inequality \eqref{applicondJensen}  then gives
\begin{multline} \label{p1terprop1induction}
 \E   \Big (\I_{ \exists i \in \{2, \dots, m_{q,p} \} \, : \, \bU_{\ell,2i}^{(k(t))} \neq  \bU_{\ell,2i}^{*(k(t))} }  {\rm Tr}\exp \big ( {\mathbb C}_{\ell-1, \ell}(t) \big )   \Big )    \\
\leq 
   \beta_{p+1} ( m_{q,p} -1)  {\rm e}^{qtM}     \,  \E   {\rm Tr}\exp \big ( {\mathbb C}_{\ell, \ell +1}(t) \big )  \, .
\end{multline} 
Starting from \eqref{p1bisprop1induction}  and taking into account \eqref{p1terprop1induction},  inequality  \eqref{p1prop1induction} follows and so does inequality \eqref{p1prop1}.

\medskip

With the same arguments as above and with obvious notations, we infer that 
\begin{multline}  \label{p2prop1} 
\E   {\rm Tr}\exp \Big ( { 2t    \sum_{j=1}^{2^{ k(t) }}  \sum_{ i=1}^{m_{q,p +1} } \bbZ_{j , 2i-1 }^{(k(t))}   } \Big ) \\
  \leq   \Big ( 1+ m_{q,p}   {\rm e}^{ 2qtM}  \beta_{p+1} \Big )^{2^{k(t)}  } 
  \E {\rm Tr}  \exp \Big (2t \sum_{j=1}^{2^{k(t)}  } \sum_{ i=1}^{m_{q,p +1} } \bbZ_{j , 2i-1 }^{*(k(t))}  \Big ) \, . 
\end{multline}
Note that to get the above inequality, we used instead of \eqref{lambdamaxpair}  that, almost surely, 
\begin{multline*} 
\lambda_{{\rm max}}  \big(2t  \sum_{ i=1}^{m_{q,p+1} } \bbZ_{\ell, 2i -1}^{(k(t))}   \big) 
\leq  2t \sum_{ i=1}^{m_{q,p +1} }  \sum_{k \in K_{k(t),\ell} \cap J^{( k(t) )}_{\ell ,2i-1}}  \lambda_{{\rm max}}\big (  \bbX_k  \big )  \\
\leq  2 M t ( m_{q,p}  p  + q - 2pm_{q,p})  = 2 M t (  q - pm_{q,p})  \leq Mt ( q +2p ) 
\leq   2 tqM \, .
\end{multline*}
Starting from \eqref{constrivial2} and taking into account \eqref{trace-exp-convex}, \eqref{p1prop1} and \eqref{p2prop1}, we then derive
\begin{multline} \label{constrivial3}
 \E {\rm Tr}\exp \Big ( t  \sum_{i \in K_{A}} \bbX_i    \Big )  
 \leq 
\Big ( 1+ m_{q,p} {\rm e}^{2qtM}  \beta_{p+1}  \Big )^{2^{k(t)}  } 
 \prod_{k=0}^{k(t) -1} \Big ( 1+ \beta_{d_k +1} {\rm e}^{tM n_{\ell} 2^{\ell - k }}\Big )^{2^k}  
  \\
 \times \Big(   \frac{1}{2} \E {\rm Tr}  \exp \big (2t \sum_{j=1}^{2^{k(t)}  } \sum_{ i=1}^{m_{q,p} } \bbZ_{j , 2i }^{*(k(t))}  \big ) 
 +   \frac{1}{2} \E {\rm Tr}  \exp \big (2t \sum_{j=1}^{2^{k(t)}  } \sum_{ i=1}^{m_{q,p} +1 } \bbZ_{j , 2i-1 }^{*(k(t))}  \big )  \Big) \, .
\end{multline}
Now we choose 
\[
p =    \Big [  \frac{2}{tM} \Big ]\wedge \Big [ \frac{q}{2}\Big ]   \, .
\]
Note that  the random vectors $(\bbZ_{j , 2i -1}^{*(k(t))} )_{i,j}$ are mutually independent and centered. Moreover, 
\[
2 \lambda_{{\rm max}} (  \bbZ_{j , 2i -1}^{*(k(t))} ) \leq   2 M p \leq  \frac{4}{t}  \ \text{ a.s.} 
\]
Therefore by using Lemma \ref{lmatropp} together with the definition of $v^2$ and the fact that $2^{ k(t)  } (m_{q,p}  +1)  p \leq 2^{ k(t)  } q  \leq A $,  we get
\beq \label{cas2e5}
\E {\rm Tr} \exp \Big({2 t  \sum_{j=1}^{2^{ k(t) }}  \sum_{ i=1}^{m_{q,p} +1} \bbZ_{j , 2i -1}^{*(k(t))}}    \Big )  
\leq  d  \exp   (  4 \times 3.1  \times A t^2 v^2  ) \, .
\eeq
Similarly, we obtain that 
\beq \label{cas2e6}
\E {\rm Tr} \exp \Big({2 t  \sum_{j=1}^{2^{ k(t) }}  \sum_{ i=1}^{m_{q,p} } \bbZ_{j , 2i }^{*(k(t))}}    \Big )  
\leq  d  \exp   (  4 \times 3.1 \times A t^2 v^2  ) \, .
\eeq
Next, by using Condition \eqref{condbeta} and \eqref{evident}, we get 
\beq \label{firstlog}
\log  \Big ( 1+ m_{q,p}   {\rm e}^{ 2tqM  }  \beta_{p+1}  \Big )^{2^{k(t)}  }  
\leq  2^{k(t) }m_{q,p}  {\rm e}^{ 2tqM  }  {\rm e}^{ -c p }   \leq \frac{A}{2p}  {\rm e}^{ 2tqM  }  {\rm e}^{ -c p }    \, .
\eeq
Several situations can occur. Either $(tM)^2 \leq \kappa/A$ and in this case $k(t) =0$ implying that $A/2 \leq q \leq A  \leq \kappa /(tM)^2$. If in addition $q \geq 4/(tM)$ then $p  = [2/(tM)]  \geq 1/tM$ (since $tM \leq 1$) and
\[
\frac{A}{2p}  {\rm e}^{ 2tqM  }  {\rm e}^{ -c p } \leq \frac{A tM }{2}    {\rm e}^{ 2 \kappa/(tM) }  {\rm e}^{ -c/(tM) }  \leq \frac{A tM }{2}   {\rm e}^{ -3c/(4tM) } \leq  \frac{ ( tM )^2 }{ c }  {\rm e}^{ -3c/(16 tM) }  \, ,
\]
where we have used that $ \log_2 A \leq \frac{ c }{32 tM} $,  $A \geq 2$, and  $ {\rm e}^{ -3c/(8tM) } \leq \frac{8tM}{3c} $ for the last inequality. 
If otherwise $q < 4/(tM)$ then $p = [q/2] \geq q/4$. Hence, since $ 2tM \leq  c/16$ (since $\log A \geq \log 2$) and $tM >4/A$, 
\[
\frac{A}{2p}  {\rm e}^{ 2tqM  }  {\rm e}^{ -c p } \leq \frac{ 2A}{q}   {\rm e}^{ -3c q/16 }\leq 4  {\rm e}^{ -3c A/32 } \leq A tM   {\rm e}^{ -3c /(8tM) }  \leq \frac{ ( tM )^2 }{  c} {\rm e}^{ -3c/(32 tM) }    \, ,
\]
where we have used that $A/2 \leq q $ for the second inequality, and that $ \log_2 A \leq \frac{ c }{32 tM} $, $A \geq 2$ and  $ {\rm e}^{ -3c/(16tM) } \leq \frac{16tM}{3c} $ for the last one. 

Either  $(tM)^2 >  \kappa/A$ and in this case $k(t) \geq 1$ and by using \eqref{evident} and the definition of $k(t)$, we have
\beq \label{consequencesurq}
q \geq \frac{A}{2^{k(t) +1} }  \geq   \frac{\kappa}{ 4 (tM)^{2}}  \, .
\eeq
If in addition $ q \geq 4/(tM)$  then $p  = [2/(tM)]  \geq 1/tM$, and   by \eqref{restnj} and the definition of $k(t)$,
\[
q \leq  2 A \frac{(1-\delta)^{\ell}}{2^{k(t)}} \leq  \frac{2\kappa}{(tM)^{2}}  \, .
\]
It follows that 
\[
\frac{A}{2p}  {\rm e}^{ 2tqM  }  {\rm e}^{ -c p }  \leq  \frac{A tM }{2}  {\rm e}^{4 \kappa/(tM)} {\rm e}^{ -c/(tM) }   \leq  \frac{A tM }{2}  {\rm e}^{ -c/(2tM) }  \leq \frac{ ( tM )^2 }{ c } {\rm e}^{ -c/(8tM) }  \, ,
\]
where we have used that $ \log_2 A \leq \frac{ c }{32 tM} $, $A \geq 2$ and  $ {\rm e}^{ -c/(4tM) } \leq \frac{4tM}{c} $ for the last inequality. 
Now if $ q <  4/(tM)$  then  $p = [q/2] \geq q/4$. Hence, using again the fact that $ 2tM \leq  c/16$ combined with \eqref{consequencesurq}, we get 
\[
\frac{A}{2p}  {\rm e}^{ 2tqM  }  {\rm e}^{ -c p } \leq \frac{ 8A (tM)^2}{\kappa}   {\rm e}^{ - 3 c q/16 }\leq \frac{ 8^2A (tM)^2}{c}   {\rm e}^{ - \frac{3 c^2 }{16 \times 4  \times 8 (tM)^2} } \leq 
\frac{ 8 (tM)^2}{  c } {\rm e}^{ -3c/(32 tM) }  \,  ,
\]
where we have used that $ \log_2 A \leq \frac{ c^2 }{(32 tM)^2} $ and $A \geq 2$ for the last inequality. 

So, overall, starting from \eqref{firstlog}, we get 
\beq \label{firstlogfinal}
\log  \Big ( 1+ m_{q,p}   {\rm e}^{ 2tqM  }  \beta_{p+1}  \Big )^{2^{k(t)}  }  
  \leq \frac{ 8 (tM)^2}{  c}  {\rm e}^{ -3c/(32 tM) }   \, .
\eeq
We handle now the term  $\prod_{k=0}^{k(t) -1} \Big ( 1+ \beta_{d_k +1} {\rm e}^{tM n_{\ell} 2^{\ell - k }}\Big )^{2^k} $ only in the case where $(\kappa/A)^{1/2} <tM$, otherwise this term is equal to one. By taking into account \eqref{condbeta}, \eqref{restdj}, \eqref{restnj} and the fact that $tM \leq  c \delta/ 8$,  we have
\begin{align*}
\log  \prod_{k=0}^{k(t) -1} \Big ( 1+ \beta_{d_k +1} {\rm e}^{tM n_{\ell} 2^{\ell - k }}\Big )^{2^k}  
& \leq    \sum_{k=0}^{k(t) -1} 2^k \exp  \Big ( -c \frac{A \delta (1 - \delta)^{k}}{2^{k+1}}  + 2 tM \frac{A (1 - \delta)^{\ell}}{2^{k}}  \Big )  
 \\&
\leq    \sum_{k=0}^{k(t) -1} 2^k \exp  \Big ( -c \frac{A \delta (1 - \delta)^{k}}{2^{k+2}}  \Big )   
\\&
\leq    2^{k(t) } \exp  \Big ( - \frac{A c \delta (1 - \delta)^{k(t) -1}}{2^{k(t)+1}}  \Big )    \, .
\end{align*}
By the definition of $k(t)$, we have $ \displaystyle 
A \frac{(1-\delta )^{k (t) - 1}}{2^{k(t) - 1}} > \frac{\kappa
}{(tM)^{2}}  $. Therefore $ 2^{k (t) } \leq   2A \frac{(tM)^{2}}{\kappa
}$.  Moreover 
\[
A  c \delta \frac{(1-\delta )^{k (t) - 1}}{2^{k (t) +1}} > \frac{ c \kappa \delta 
}{4 (tM)^{2}} \geqslant \frac{ 2 \kappa 
}{tM} \, ,
\]
since $tM \leq  c \delta/ 8$.  It follows that 
\beq \label{cas2e8}
\log  \prod_{k=0}^{k(t) -1} \Big ( 1+ \beta_{d_k +1} {\rm e}^{tM n_{\ell} 2^{\ell - k }}\Big )^{2^k}  \leq      2A \frac{(tM)^{2}}{\kappa} \exp  \big ( - 2 \kappa /(tM) 
\big )   \leq         \frac{  (tM)^{2}}{ c } {\rm e}^{ -3c/(32 tM) } 
  \, ,
\eeq
where we have used the fact that $ \log_2 A \leq \frac{ c }{32 tM} $.  So, overall, starting from \eqref{constrivial3} and considering the upper bounds \eqref{cas2e5}, \eqref{cas2e6}, \eqref{firstlogfinal} and  \eqref{cas2e8}, we get 
\[
\log \E {\rm Tr}\exp \Big ( t  \sum_{i \in K_{A}} \bbX_i    \Big )  
 \leq  \log d    +  4 \times 3.1 A t^2 v^2   +   \frac{9 (tM)^2}{ c }   {\rm e}^{ -3c/(32 tM) }   \, .
\]
Therefore Inequality \eqref{resultpropinter1} also holds in the case where   $4/A <t M \leq  \min  \big ( \frac{1}{2}, \frac{c \log 2}{32 \log A} \big ) $.   This ends the proof of the proposition. $\square$

\subsection{Proof of Theorem \ref{mainthm}}  \label{sectionend}

 Let $A_{0}=A=n$ and $%
\bbY^{(0)}(k)=\bbX_{k}$ for any $k=1,\dots ,A_{0}$. Let $K_{A_{0}}$ be the discrete Cantor type set as defined from $%
\{1,\dots ,A\}$ in Section \ref{sectionCantor}. Let $%
A_{1}=A_{0}- {\rm Card} (K_{A_{0}})$ and define for any $k=1,\dots
,A_{1}$, 
\begin{equation*}
\bbY^{(1)}(k)=\bbX_{i_{k}}\text{ where }\{i_{1},\dots ,i_{A_{1}}\}=\{1,\dots
,A\}\setminus K_{A}\,.
\end{equation*}%
Now for $i\geq 1$, let $K_{A_{i}}$ be defined from $\{1,\dots
,A_{i}\}$ exactly as $K_{A}$ is defined from $\{1,\dots ,A\}$.
Set $A_{i+1}=A_{i}-{\rm Card} (K_{A_{i}}) $ and $\{j_{1},\dots
,j_{A_{i+1}}\}=\{1,\dots ,A_{i}\}\setminus K_{A_{i}}$. Define now 
\begin{equation*}
\bbY^{(i+1)}(k)=\bbY^{(i)}(j_{k})\text{ for }k=1,\dots ,A_{i+1}\,, 
\end{equation*}%
and set
$$
L=L_n = \inf \{ j \in {\mathbb N}^* \, , \, A_j \leq  2 \}
\, .
$$
Note that, for any $i \in \{0, \dots,  L-1\}$,  $A_i >2$ and  ${\rm Card } (K_{A_i}) \geqslant A_i/2$. Moreover $A_i \leq  n 2^{-i}$. 

The following decomposition clearly 
holds
\begin{equation} \label{P3prop3}
\sum_{k=1}^n \bbX_{k} = \sum_{i=0}^{L-1} \sum_{k \in K_{A_i}} \bbY^{(i)} (k)  +
\sum_{k=1}^{A_L} \bbY^{(L)} (k) \, .
\end{equation}
Let
\[
{\mathbb U}_i= \sum_{k \in K_{A_i}} \bbY^{(i)} (k)  \text{ for $0\leq i \leq  L-1$ and
}{\mathbb U}_L=\sum_{k=1}^{A_L} \bbY^{(L)} (k)  \, ,
\]
For any positive $x$,  let 
 \[
 h(c, x) = \min  \Big ( \frac{1}{2}, \frac{c \log 2}{32 \log x} \Big ) \, .
 \]
For any $i \in \{0, \dots,  L-1\}$, noticing that the self-adjoint random matrices $(Y^{(i)} (k) )_k$ satisfy the condition \eqref{condbeta} with the same constant $c$, we can apply Proposition \ref{propinter1} and  get that  for any positive $t$ satisfying $t M< h (c,n/ 2^{i} )$,
\beq  \label{appli1}
\log \E {\rm Tr}\big ( \exp (t  {\mathbb U}_i  ) \big )    \leq  \log d +  \frac{ 4 t^2 n 2^{-i}( 2 v + \sqrt{3} \times 2^{5i/2}M / (n^{5/2} {\sqrt c}))^2 }{1 - t M /h (c,n 2^{-i} ) } \, .
\eeq
On the other hand, by Weyl's inequality,
\[
\lambda_{{\rm max}} \big (  {\mathbb U}_L \big ) \leq  M A_L \leq  2 M  \, .
\]
Therefore by using Lemma \ref{lmatropp}, for any positive $t$,
\[
 \E {\rm Tr}\big ( \exp (t  {\mathbb U}_L  ) \big )   \leq  d  \exp  \Big (    t^2 g(2tM)  \lambda_{{\rm max}}  \big (\E \big (   {\mathbb U}_L^2 \Big )  \big )\, .
\]
Hence by the definition of $v^2$, for any positive $t$ such that $tM < 1$, we get 
\beq  \label{appli2}
\log  \E {\rm Tr}\big ( \exp (t  {\mathbb U}_L  ) \big )  \leq  \log d +     2 t^2  v^2  \leq   \log d +  \frac{ 2 t^2  v^2 }{1 -  t M   } \, .
\eeq
Let 
\[
{ \kappa}_i=  \frac{M}{ h (c,n/ 2^{i} ) } \text{ for $0\leq i \leq  L-1$ and
}{\kappa}_L=M    
\]
and 
\[
{ \sigma}_i=  2 \frac{\sqrt{ n}}{ 2^{ i/2}} \Big (  2v +  \sqrt{3} \times \frac{2^{i}M}{n {\sqrt c}} \Big )  \text{ for $0\leq i \leq  L-1$ and
}{\sigma}_L= v \sqrt{2}  \, .
\]
Since 
\[
L \leq  \Big [ \frac{ \log n - \log 2}{\log 2}\Big ] +1 \, , 
\]
we get 
\[
\sum_{i=0}^L { \kappa}_i \leq  M    \Big (  \sum_{i=0}^{L-1} \frac{1}{h (c,n/ 2^{i} ) } + 1  \Big ) 
 \leq   M \frac{\log n}{\log 2} \max \Big ( 2 , \frac{32 \log n}{c \log 2}  \Big )  = M  \gamma (c, n)   \, .
\]
Moreover
\begin{align*}
\sum_{i=0}^L { \sigma}_i & = 2 \sqrt{ n} \sum_{i=0}^{L-1} 2^{-i/2}\Big (  2 v + \sqrt{3} \times \frac{2^{5i/2}M}{n^{5/2} {\sqrt c}} \Big )  \,  + \,   v \sqrt{2}  \\
& \leq   14 \sqrt{ n }  v  +    2 c^{-1/2}n^{-2}   M 2^{2L} +  v \sqrt{2}    \\
& \leq    15 \sqrt{ n }  v  +     2 c^{-1/2}   M   \, .
\end{align*}
Taking into account \eqref{appli1} and \eqref{appli2}, we get overall by Lemma \ref{breta}, that for any positive $t$ such that $\displaystyle tM < 1/\gamma (c, n) $,
\beq \label{resloglaplace}
\log \E   {\rm Tr} \Big (\exp \big ( t \sum_{i =1}^n\bbX_i    \big ) \Big ) \leq  \log d +  \frac{     t^2 n  \big ( 15 v + 2 M/ (cn)^{1/2} \big )^2 }{1 - t M \gamma (c, n)  } := \gamma_n(t)\, .
\eeq
To end the proof of the theorem, it suffices to notice that for any positive $x$
\[
{\mathbb P} \Big (  \lambda_{{\rm max}} \big ( \sum_{i =1 }^n \bbX_i  \big ) \geqslant x  \Big )  \leq  \inf_{t  >0 \, : \, tM \leq  1/\gamma (c, n)  } \exp \big ( - t x + \gamma_n(t) \big ) \, ,
\]
where $\gamma_n(t) $ is defined in \eqref{resloglaplace}. $\square$

\medskip

\noindent{\bf Acknowledgements.} This work was initiated when the first named author was visiting the third named author at the University of Alberta. She would like to thank Nicole Tomczak-Jaegermann
and Alexander Litvak for their hospitality and the University of Alberta for the excellent working conditions.


\medskip
\begin{thebibliography}{99}

\bibitem{Adam}  Adamczak, R. (2008). A tail inequality for suprema of unbounded empirical processes with applications to Markov chains. \textit{Electron. Journal Probab.}\textbf{ 13},  1000-1034.

\bibitem{AhWi} Ahlswede, R. and   Winter, A. (2002). Strong converse for identification via quantum channels. \textit{IEEE Trans. Inform. Theory}, \textbf{48}, no. 3, 569-579. 

\bibitem{Be} Berbee H.C.P. (1979). Random walks with stationary increments and renewal theory. \textit{%
Mathematical Centre Tracts},  \textbf{112},  Mathematisch Centrum,  Amsterdam, 223 pp.


\bibitem{Ch-thesis} Chatterjee, S. (2006). Concentration inequalities with exchangeable pairs. Ph.D. thesis, \textit{Stanford Univ., Palo Alto.}.

\bibitem{Ch} Chatterjee, S. (2007). Stein’s method for concentration inequalities. \textit{Probability theory and related fields}, \textbf{138}, no. 1, 305-321.

\bibitem{Do} Doob, J. L. (1953) {\it  Stochastic processes}. Wiley, New York.

\bibitem{Ibra}  Ibragimov, I. A. (1962). Some limit theorems for stationary processes. \textit{Teor. Verojatnost. i Primenen}. {\bf 7},  361-392.



\bibitem{KR} Kolmogorov, A.N., Rozanov, Y.A. (1960). On the strong mixing conditions for stationary
gaussian sequences.  {\it Theor. Probab. Appl.}, {\bf 5}, 204-207.



\bibitem{Li} Lieb, E. H. (1973). Convex trace functions and the Wigner-Yanase-Dyson conjecture. \textit{Adv. Math.} \textbf{11}, 267-288.


 

\bibitem{Ma et al} Mackey, L., Jordan, M. I., Chen, R. Y., Farrell, B.,  Tropp, J. A. (2014). Matrix concentration inequalities via the method of exchangeable pairs. \textit{ Ann. Probab.}, \textbf{42}, no. 3, 906-945.


\bibitem{MPR1}  Merlev\`ede F., Peligrad M. and Rio E. (2009). Bernstein inequality and moderate deviations under strong mixing conditions. \textit{High dimensional probability V: the Luminy volume, 273–292, Inst. Math. Stat. Collect., 5, Inst. Math. Statist., Beachwood, OH}
 
 \bibitem{MPR2}  Merlev\`ede F., Peligrad M. and Rio E. (2011). A Bernstein type inequality and moderate deviations for weakly dependent sequences. \textit{Probab. Theory Related Fields} \textbf{151}, no. 3-4, 435-474. 
 
 \bibitem{Mo} Mokkadem, A. (1990). Propri\' et\' e􏰌 de m\' e􏰌lange des processus autor\' egressifs polynomiaux. {\it Ann.
Instit Henri Poincar\'e.} {\bf 26}, 133-141.
 

\bibitem{Pa et al} Paulin, D., Mackey, L. and   Tropp, J. A. (2014). Efron-Stein Inequalities for Random Matrices. \textit{arXiv:1408.3470}.



\bibitem{Sk} Skorohod, A. V.  (1976). On a representation of random variables. \textit{Theory Probab. Appl.}, \textbf{21}, 628-632.

\bibitem{Tao12} Tao, T. (2012). {\it Topics in random matrix theory}. Graduate Studies in Mathematics, 132. American Mathematical Society, Providence, RI, 282 pp.


\bibitem{Tr} Tropp, J. A. (2012). User-friendly tail bounds for sums of random matrices. \textit{ Found. Comput. Math.}, {\bf 12}, no. 4, 389-434.



\end{thebibliography}
\end{document}